\documentclass[11pt,a4paper]{article}
\usepackage{amssymb,latexsym,a4wide}
\def\ZZ{{\mathbb Z}}
\def\QQ{{\mathbb Q}}
\def\CC{{\mathbb C}}
\def\proof{{\bf Proof}.\ }
\def\bull{\vrule height .9ex width .8ex depth -.1ex }

\newtheorem{formula}{}[section]
\newtheorem{proposition}[formula]{Proposition}
\newtheorem{definition}[formula]{Definition}
\newtheorem{corollary}[formula]{Corollary}
\newtheorem{remark}[formula]{Remark}
\newtheorem{lemma}[formula]{Lemma}
\newtheorem{theorem}[formula]{Theorem}

\def\thrm{\begin{theorem}}
\def\thrml#1{\begin{theorem}\label{#1}}
\def\ethrm{\end{theorem}}
\def\rmrk{\begin{remark}}
\def\rmrkl#1{\begin{remark}\label{#1}}
\def\ermrk{\end{remark}}
\def\dfntn{\begin{definition}}
\def\dfntnl#1{\begin{definition}\label{#1}}
\def\edfntn{\end{definition}}
\def\nmrt{\begin{enumerate}}
\def\enmrt{\end{enumerate}}

\def\qtn{\begin{equation}}
\def\qtnl#1{\begin{equation}\label{#1}}
\def\eqtn{\end{equation}}
\def\lmm{\begin{lemma}}
\def\lmml#1{\begin{lemma}\label{#1}}
\def\elmm{\end{lemma}}
\def\crllr{\begin{corollary}}
\def\crllrl#1{\begin{corollary}\label{#1}}
\def\ecrllr{\end{corollary}}
\begin{document}
\title{Analogue of Newton-Puiseux series for non-holonomic D-modules and factoring}
\author{
Dima Grigoriev \\[-1pt]
\small CNRS, Math\'ematiques, Universit\'e de Lille \\[-3pt]
\small 59655, Villeneuve d'Ascq, France\\[-3pt]
{\tt \small dima@mpim-bonn.mpg.de}\\[-3pt]
\small http://logic.pdmi.ras.ru/\~{}grigorev
}
\date{}
\maketitle
\begin{abstract}
We introduce a concept of a fractional-derivatives series and
prove that any linear partial differential equation in two
independent variables has a fractional-derivatives series solution
with coefficients from a differentially closed field of zero
characteristic. The obtained results are extended from a single
equation to $D$-modules having infinite-dimensional space of
solutions (i.~e. non-holonomic $D$-modules). As applications we design algorithms for treating first-order
factors of a linear partial differential operator, in particular for finding all
(right or left) first-order factors.
\end{abstract}

\vspace{4mm}
{\bf keywords:} Newton-Puiseux series for $D$-modules, fractional derivatives,
factoring linear partial differential operators

\vspace{4mm}

\section*{Introduction}
It is well-known that any polynomial equation $t(x,y)=0$ has
$deg_y(t)$ (counting with multiplicities) zeroes being
Newton-Puiseux series (see e.~g. \cite{Walker})
\begin{equation}\label{0}
y(x)=\sum_{i_0\leq i < \infty} y_ix^{-i/q}
\end{equation}
\noindent
for suitable integers $q\geq 1,i_0$ and the coefficients $y_i$
from an algebraically closed field.

In this paper an analogue of Newton-Puiseux series for partial
linear differential equations $T=0$ is proposed, and we prove that
$T=0$ has a solution of this form. Whereas a Newton-Puiseux series
is developed for a (plane) curve, we restrict ourselves with linear
partial differential operators $T$ in two derivatives $d_x,d_y$ (in
case of 3 or more derivatives there are no solutions of this form in
general, see Remark~\ref{three}).

One of the principal features of Newton-Puiseux series is the
appearance of fractional exponents. Thus, a question arises, what
could be an analogue of fractional powers, so to say "fractional
derivatives"? An evident observation shows that in the derivative
$y'(x)=\sum_i (-i/q+1)y_{i-q}x^{-i/q}$ the $i$-th coefficient depends
on the $(i-q)$-th coefficient of $y(x)$ itself.

That is why as a differential analogue of Newton-Puiseux series we
suggest a {\it fractional-derivatives series} of the form
$$\sum_{0\leq i<\infty} h_iG^{(-i/q)}$$
\noindent
where $h_i$ being elements of a differentially closed (or
universal in terms of \cite{Kolchin}) field $F$ and $G^{(-i/q)}$ is
called {\it $(-i/q)$-th fractional derivative} of $G$. The symbol
$G=G^{(0)}=G_{(s_2,\dots,s_k)} (f_1,f_2,\dots,f_k)$ is defined by
rational numbers $1>s_2>\cdots>s_k>0$ and $f_1,\dots,f_k\in F$
(if to continue the analogy with curves, $G$ plays a role of a uniformizing element).
For any rational $s$ the $s$-th fractional derivative $G^{(s)}$ fulfills the
identity
$$dG^{(s)}=(df_1)G^{(1+s)}+(df_2)G^{(s_2+s)}+\cdots+(df_k)G^{(s_k+s)}$$
\noindent
where either a derivative $d=d_x$ or $d=d_y$. The common
denominator $q$ of $s_2,\dots,s_k$ plays a role similar to one of
the common denominator of the exponents in a Newton-Puiseux series
(\ref{0}). The inequality $k\leq q$ holds.

In a particular case $k=1$ we have $q=1$ and as $G$ one can take $g(f_1)$
for any univariate ("undetermined") function $g$, provided that the
composition makes sense, the fractional derivatives $G^{(s)}=g^{(s)}(f_1)$ for
integers $s$. We note that finite sums
$$\sum_{i_0\leq i\leq i_1} h_iG^{(-i)}$$
\noindent
(so, for $k=q=1$) appear in the Laplace method as solutions of some
second-order equations $T=0$ (see e.~g.
\cite{Goursat, Tsarev}).

One can find necessary in the sequel information on $D$-modules in
\cite{Bjork}, \cite{Sabbah}, a survey on their
algorithmical aspects in \cite{Quadrat}. We mention that there are also
applications of Newton polygons over the Weyl algebra $C[x,d_x]$: in
\cite{Sabbah} to meromorphic connections, in \cite{Malgrange81} to
micro-differential operators and in \cite{Malgrange91} to the Fourier
transform.
In case of linear {\it ordinary} differential operators Newton polygons are
employed to produce the canonical form basis of the space of solutions (see
e.~g. \cite{Wasow}, also \cite{G90} where an algorithm for this problem with
a better complexity bound was designed). A similar form of solutions for
linear  {\it partial} differential operators were studied in \cite{Cano} where,
nevertheless, also examples are exhibited of operators without solutions of
this form.
On the problem of factoring a linear ordinary
differential operator one can look in \cite{Singer2003}, see also \cite{G90}.

In Section~\ref{series} we introduce the principal concept of
fractional-derivatives series and give some their basic
properties.

In the sequel the crucial role plays the multiplicity $m$ of a
linear factor of the {\it symbol} of the linear partial differential
operator $T$ (with coefficients in $F$) of an order $n$: the
symbol is a homogeneous polynomial in two variables $d_xf_1,d_yf_1$ of the degree
$n$ which corresponds to the highest derivatives of $T$. In
Section~\ref{solution} we develop a method for constructing
fractional-derivatives solutions of $T=0$ and prove the existence
of such a solution with $q\leq m$. The method is similar to the
Newton-Puiseux expansion, it produces a relevant convex polygon
similar to the Newton one, but differs in several aspects. The
main of the latter is that the leading equation corresponding to
a certain (leading) edge of the polygon is not a univariate polynomial
unlike the Newton-Puiseux expansion, but rather a non-linear
first-order partial differential equation. This creates
difficulties in defining a multiplicity of a solution of the
leading equation. Also it is unclear, what could be a
differential analogue of the statement (cf. above) that
an algebraic equation $t=0$ has precisely $deg_y(t)$
Newton-Puiseux series
solutions of the form (\ref{0})? Partially these questions are answered for
the introduced in Section~\ref{generic} generic fractional-derivatives series
solutions.

In Section~\ref{ideal} the result of Section~\ref{solution} is
extended from a single partial linear differential equation to a
system of equations in several unknown functions having an
infinite-dimensional space of solutions (or in other words, to a
$D$-module of a non-zero differential type, one can call it a
non-holonomic $D$-module). To this end for any left ideal $J\subset
F[d_x,d_y]$ of the differential type 1 we yield an operator $p\in
F[d_x,d_y]$ and show that any fractional-derivatives series solution
of the equation $p=0$ which corresponds to a linear factor
(different from $d_yf_1$) of the symbol of $p$, is a solution of the
ideal $J$ as well. In Section~\ref{duality} we exploit the relation
of equivalence of ideals introduced in \cite{GS05} and establish a
kind of duality between equivalence classes of non-holonomic ideals
and their sets of fractional-derivatives series solutions. Namely,
it is proved that two non-holonomic left ideals $J, J_1 \subset
F[d_x,d_y]$ are equivalent if and only if their respective sets of
fractional-derivatives series solutions coincide. Also we express
the quotient of the spaces of fractional-derivatives series
solutions of non-holonomic ideals $J\subset J_1$ via the module of
relative syzygies \cite{GS05} of this pair of ideals.

In Section~\ref{completeness} it is shown that in case of a separable operator $T$
any its power series solution can be obtained as a sum of specifications of its
suitable fractional-derivatives series solutions, thereby establishing completeness
of the latter. In Section~\ref{applications} we provide applications of fractional-derivatives series
to studying first-order factors of an operator, exploiting that in case of a first-order
operator $T=d_y+ad_x+b$ its fractional-derivatives series solutions turn to a single term
of the form $hG(f)$ where $(d_y+ad_x)f=0$ and $h$ being a particular solution of $T=0$.
In Subsection~\ref{factors} an algorithm is designed which finds first-order factors of a
given operator, and in Subsection~\ref{radical} an algorithm which constructs the intersection
of all the principal ideals generated by the first-order factors of the operator. In Section~\ref{second} the possible fractional-derivatives series solutions
of a second-order operator obtained by the algorithm from Section~\ref{solution} are described. This description can help to imagine the shape of
fractional-derivatives series solutions and the difficulties which appear
while their developing.

\section{Fractional-derivatives series}\label{series}
Let $F$ be a differential field of the characteristic 0 with the
derivatives $\{d_j\}$ and
a subfield of constants $C\subset F$ \cite{Kolchin}.

\begin{definition}\label{symbol}
Let $f_1,\dots,f_{k_0} \in F$ and rational numbers $1 > s_2 >
\cdots s_{k_0} > 0$. We introduce a symbol
$G=G^{(0)}=G_{s_2,\dots,s_{k_0}}(f_1,f_2,\dots,f_{k_0})$ together a
set $\{G^{(s)}\}_{s\in \QQ}$ of its fractional $s$-th derivatives satisfying
the following
rule of differentiation for any derivative $d=d_j$:
$$dG^{(s)}=(df_1)G^{(1+s)}+(df_2)G^{(s_2+s)}+\cdots+(df_{k_0})G^{(s_{k_0}+s)}$$
\end{definition}

Clearly, these differentiations commute with each other and one
can consider the free $F$-module with the basis $\{G^{(s)}\}_{s\in
\QQ}$ as a $D$-module.

\begin{definition}
Let $q$ be the common denominator of $s_2,\dots,s_{k_0}$ and $h_i\in
F, i\geq 0, s_0q\in \ZZ$. Then
\qtnl{1}
H=\sum_{0\leq i < \infty} h_i G^{(s_0-i/q)}
\eqtn
\noindent
we call a {\it fractional-derivatives series}.
\end{definition}

For a given $G$ all the fractional-derivatives series (with added 0) constitute a
$D$-module (we study it below in Section~\ref{ideal}). Obviously, $k_0\leq q$.

It is easy to see that $G$ satisfies a suitable linear partial
differential equation with coefficients in $F$.

\begin{remark}\label{specialisation}
The symbol $G$ plays a role in $H$ similar to the role of the
parameter $x$ in a Newton-Puiseux series (\ref{0}). In particular,
specifying the values of $x$ in a certain field one gets points of
(a branch of) the curve given by (\ref{0}). Here one can also
provide some specifications of $G$. Indeed, for an arbitrary
family $\{c_{i/q}\}_{i\in \ZZ}$ where $c_{i/q}\in \CC$ the
following set
$$G^{(s)}=\sum_{j_1\geq 0,\dots, j_{k_0}\geq 0}
c_{-s-j_1-j_2s_2-\cdots -j_{k_0}s_{k_0}} {f_1^{j_1} \over j_1!}
\cdots {f_{k_0}^{j_{k_0}} \over j_{k_0}!}$$
\noindent
satisfies Definition~\ref{symbol}.

For example, in case when $F$ is the ring of functions analytic in
a certain neighborhood of a given point in the multidimensional complex
space and the absolute values $|c_{i/q}|$ are bounded, the
latter series also converges in a suitable neighborhood.
\end{remark}

From now on let $F$ have two derivatives $\{d_x,d_y\}$. Consider a
linear operator

\qtnl{1.5}
T=T_0+\cdots+T_n
\eqtn
\noindent
of the order $n$ where $T_p=\sum_{0\leq j\leq p}
b_{j,p}d_x^jd_y^{p-j}$ contains the derivatives of the order $p$
and the coefficients $b_{j,p}\in F$. The following lemma holds, in
fact, for an arbitrary number of derivatives, nevertheless, the
assumption that $F$ has  two derivatives simplifies the
notations and in the sequel we deal just with operators in two
derivatives (one can verify lemma by a direct calculation).

\begin{lemma}\label{sum}
$d_x^jd_y^{p-j}(hG)$ equals the sum of the terms of the form

$${1 \over w_1!\cdots w_N!} {j \choose
l_{1,1},\dots,l_{1,m_1},\dots,l_{k_0,1},\dots,l_{k_0,m_{k_0}},l}
{p-j \choose
r_{1,1},\dots,r_{1,m_1},\dots,r_{k_0,1},\dots,r_{k_0,m_{k_0}},r}$$

$$\prod_{1\leq i\leq m_1} (d_x^{l_{1,i}}d_y^{r_{1,i}}f_1)\cdots
\prod_{1\leq i\leq m_{k_0}} (d_x^{l_{k_0,i}}d_y^{r_{k_0,i}}f_{k_0})
(d_x^ld_y^rh) \cdot G^{(m_1+s_2m_2+\cdots+s_{k_0}m_{k_0})}$$

\noindent for all partitions
$l_{1,1}+\cdots+l_{1,m_1}+\cdots+l_{k_0,1}+\cdots+l_{k_0,m_{k_0}}+l=j$
of $j$ and
$r_{1,1}+\cdots+r_{1,m_1}+\cdots+r_{k_0,1}+\cdots+r_{k_0,m_{k_0}}+r=p-j$
of $p-j$ such that $l_{\kappa, i}+r_{\kappa, i}\ge 1$ for every
$1\le \kappa \le k_0, \, 1\le i\le m_{\kappa}$ , where $w_1,\dots,
w_N$ denote the cardinalities of the partition of the triples
$(l_{1,1}, \, r_{1,1}, \, 1),\dots , (l_{1,m_1}, \, r_{1,m_1}, \,
1), \dots , (l_{k_0,1}, \, r_{k_0,1}, \, k_0), \dots ,
(l_{k_0,m_{k_0}}, \, r_{k_0,m_{k_0}}, \, k_0)$
into equal
ones, in particular, $w_1+\cdots+w_N=m_1+\cdots+m_{k_0}$.
\end{lemma}

\section{Constructing fractional-derivatives series
solutions}\label{solution}

From now on we suppose that the field $F$ is differentially closed
(or universal in terms of
\cite{Kolchin}).

The main purpose of this section is to prove that a linear partial
differential equation $T=0$, see (\ref{1.5}), has a solution of the
form (\ref{1}). To simplify the notations we put $s_0=0$ and
$h=h_0\neq 0$ in (\ref{1}).

Denote by ${\bar T_p}(d_xf_1,d_yf_2)=\sum_{0\leq j\leq
p}b_{j,p}(d_xf_1)^j(d_yf_1)^{p-j}$ a homogeneous form of the
degree $p$ in $d_xf_1,d_yf_1$. Sometimes, ${{\bar T_n}=symb(T)}$ is called the
symbol of $T$. Fix a linear factor $a_1d_xf_1+a_2d_yf_1$ of ${\bar T_n}$
having a multiplicity $m$, the coefficients $a_1,a_2\in F$.

Expanding $T(H)$ with respect to the fractional derivatives $\{G^{(s)}\}_s$
for $k=1$ (in
other words, assuming for the time being that
$dG^{(s)}=(df_1)G^{(1+s)}$, see Definition~\ref{symbol}), we get
that the coefficient at $G^{(n)}$ vanishes, i.~e. $h\cdot symb(T)=0$.
Thus, we can suppose that
$(a_1d_x+a_2d_y)f_1=0$. Choose any such $f_1$ with $grad(f_1)\neq
0$.

For $k\geq 2$ we introduce an auxiliary polygon $P_k$ playing the
role similar to the Newton polygon. Now let $k=2$, in other words, we assume
(for the time being) that
$dG^{(s)}=(df_1)G^{(1+s)}+(df_2)G^{(s_2+s)}$. The next purpose is
to construct $s_2$ and $f_2$. It suffices to consider the
expansion of the first term $T(hG)$ of $T(H)$
(we'll come back to this issue at the end of
the present section). When we talk about the expansion of $T(hG)$ we always refer
to Lemma~\ref{sum}. If a term $b(\prod_{1\leq i\leq t}
(d_x^{l_i}d_y^{r_i}f_2))G^{(s+s_2t)}$ occurs in $T(hG)$, where $b$
is a differential polynomial in $f_1$ and in $h$ (being linear in
$h$), then we place the point $(s,t)$ in $P_2$. As $P_2$ we take
the convex hull of these points with the origin $(0,0)$. If to
assign the weight 1 to every derivative $d_x^ld_y^rf_1$ then any
term in $b$ gets the weight $s$ due to Lemma~\ref{sum}.

One can observe that $P_2$ lies to the left from the line ${\bar
L_1}=\{s+t=n\}$ with the slope 1 (under the slope of the line
$\{s+jt=const\}$ we mean $j$) again due to  Lemma~\ref{sum}. Moreover,
the point $(n-m,m)\in {\bar L_1}$ belongs to $P_2$ because the
non-zero term
$${{\bar T_n} \over ((a_1d_x+a_2d_y)f_1)^m} \cdot ((a_1d_x+a_2d_y)f_2)^m
\cdot G^{(n-m+s_2m)}$$

\noindent
occurs in the expansion of $T(hG)$, taking into account that the factor
$(a_1d_x+a_2d_y)f_1$ has the multiplicity $m$ in ${\bar T_n}$, and
no other term from this expansion gives a contribution in the
coefficient at the point $(n-m,m)$. Similarly, one verifies that the points
$(n-t,t)$ with $0\leq t\leq m-1$ do not belong to $P_2$.

Now we assign a (yet unknown) weight $s_2$ to every derivative
$d_x^ld_y^rf_2$. Therefore, to find $s_2<1$ we consider the edges
of $P_2$ with the positive slopes less than 1. Choose any such
edge $L_2$ (we call it {\it leading}) with the endpoints
$(j_1,t_1), \, (j_2,t_2), \, t_1>t_2$; we have seen already that
$t_1\leq m$. Then the slope of $L_2$ provides
$s_2=(j_1-j_2)/(t_1-t_2)$.

To find $f_2$ we consider the {\it leading} differential
polynomial $Q_2(f_2)$ which equals the sum of the coefficients at
all the points of $P_2$ which lie on $L_2$. Then  $Q_2(f_2)$
coincides with the coefficient at $G^{(j_1+s_2t_1)}$ in the
expansion of $T(hG)$. As $f_2\in F$ we take a solution of the
leading equation $Q_2(f_2)=0$. Evidently, $j_1+s_2t_1< n$
since the point of intersection of the line $\bar L_2$ (which
contains the edge $L_2$) with $j$-axis $\{t=0\}$ is located to the
left of the intersection of $\bar L_1$ with $j$-axis.

Thus, we are able to formulate the recursive hypothesis of the
procedure under description which constructs $1>s_2>s_3>\cdots$
and $f_1,f_2,f_3,\dots$. Suppose that $s_2,\dots,s_k$ and
$f_1,f_2,\dots,f_k$ are already constructed. In addition, a polygon
$P_k$ is constructed being a convex hull of the points $(j,t)$
(together with the origin $(0,0)$) such that a term
\qtnl{2}
b(\prod_{1\leq i\leq t} d_x^{l_i}d_y^{r_i}f_k)  G^{(j+s_kt)}
\eqtn
\noindent
occurs in the expansion of $T(hG)$ under the assumption
$dG=(df_1)G^{(1)}+(df_2)G^{(s_2)}+\cdots+(df_k)G^{(s_k)}$, see
Definition~\ref{symbol}. A
certain leading edge $L_k$ of $P_k$ is chosen with a slope $s_k>0$
and with the endpoints $(j_3,t_3), \, (j_4,t_4), \, t_3>t_4$. We
name $(j_3,t_3)$ the {\it pivot} of $L_k$ and $t_3$ the {\it
multiplicity} of $L_k$. The leading differential polynomial
$Q_k(f_k)$ equals the sum of the coefficients at all the points of
$P_k$ which lie on $L_k$. Then $Q_k(f_k)$ coincides with the
coefficient at $G^{(j_3+s_kt_3)}$ in the expansion of $T(hG)$. As
$f_k\in F$ a solution of the leading equation  $Q_k(f_k)=0$
is taken.
The points of intersections of the lines ${\bar L_1}, {\bar L_2},
\dots$ with $j$-axis decrease. Denote by $q_k$ the common
denominator of $s_2,\dots, s_k$, obviously $q_1=1$.

To carry out the recursive step, we make the assumption
$dG=(df_1)G^{(1)}+(df_2)G^{(s_2)}+\cdots+(df_k)G^{(s_k)}+
(df_{k+1})G^{(s_{k+1})}$. The boundary of the polygon $P_{k+1}$
above the pivot of $L_k$ (including the pivot itself) is the same
as of $P_k$.

Let us calculate the points of $P_{k+1}$ located on the line $\bar
L_k$. Denote by $B_t \, (t_4\leq t\leq t_3)$ the coefficient of
$P_k$ at the point $(j_3+s_k(t_3-t),t)\in L_k$. Then
$Q_k=\sum_{t_3\leq t\leq t_4} B_t$. One can observe that $B_t$
contains no higher derivative $d_x^ld_y^rf_k$ with $l+r\geq 2$.
Indeed, if otherwise $B_t$ contained a term of the form (\ref{2})
then the coefficient of $P_k$ at the point
$(j_3+s_k(t_3-t),t+\sum_{1\leq i\leq t} (l_i+r_i-1))$ would
contain the term
$$b(d_xf_k)^{\sum_{1\leq i\leq t} l_i}(d_yf_k)^{\sum_{1\leq i\leq
t} r_i} G^{(j_3+s_k(t_3+\sum_{1\leq i\leq t} (l_i+r_i-1)))}$$
\noindent
due to Lemma~\ref{sum}, hence the point $(j_3+s_k(t_3-t),
t+\sum_{1\leq i\leq t} (l_i+r_i-1))$ should belong to $P_k$ which
leads to a contradiction when $\sum_{1\leq i\leq t} (l_i+r_i-1)
\geq 1$.

Besides, $B_t$ is a linear form in the derivatives of
$h$. We claim that $B_t=h{\tilde B_t}$ for an appropriate
differential polynomial $\tilde B_t$ in $f_1, \dots , f_k$. Indeed, if $B_t$
contained a term $(d_x^ld_y^rh){\tilde b}G^{(j_3+s_kt_3)}$ with
$l+r\geq 1$ for a certain $\tilde b$ being a differential
polynomial in $f_1,\dots, f_k$, then the coefficient of $P_k$ at
the point $(j_3+s_k(t_3-t),t+l+r)$ would contain the term
$h{\tilde b}(d_xf_k)^l(d_yf_k)^rG^{(j_3+s_k(t_3+l+r))}$ due to Lemma~\ref{sum},
therefore, the point $(j_3+s_k(t_3-t),t+l+r)$ should belong to $P_k$, the
achieved contradiction proves the claim.

Thus, $B_t$ will be treated as a homogeneous (of the degree $t$)
polynomial in $d_xf_k,d_yf_k$. For more generality of the
auxiliary results below we deem that $B_t$ is a homogeneous
polynomial in the variables $v_1,\dots, v_p$, thereby $p=2$ and
$v_1=d_xf_k, \, v_2=d_yf_k$. We denote the corresponding
derivatives ${\bar v_1}=d_xf_{k+1}, \, {\bar v_2}=d_yf_{k+1}$.

\begin{remark}\label{canonical}
Since the main purpose of the present section is to prove the existence
of solutions of the form (\ref{1}) of an equation $T=0$ (see
(\ref{1.5})) it suffices to study only the {\it canonical} solutions, namely,
when each $s_k$ is the slope of a certain edge of $P_k$ and $f_k$
satisfies a leading equation. Alternatively, one could take
$s_k$ to be the slope of some line passing through a single
vertex, say $(j_3,t_3)$ of $P_k$. In this case $B_{t_3}(f_k)=0$,
because $B_{t_3}$ is a homogeneous polynomial in $d_xf_k, d_yf_k$,
we get that $f_k$ fulfills a certain first-order linear equation
$b_1d_xf_k+b_2d_yf_k=0$. There is no way to bound the denominators
$s_k$ for non-canonical solutions (\ref{1}), the number of steps
$k_0$, moreover, the procedure of constructing $1>s_2>s_3>\cdots$
and $f_1,f_2,f_3, \dots$ could last infinitely. One might even choose
real exponents $s_k$ (cf. \cite{GS91} where an analogue of
Newton-Puiseux series solutions with real exponents was studied
for {\it non-linear} ordinary differential equations).
\end{remark}

Denote by ${\bar B_t}(0\leq t\leq t_3)$ the coefficient at the
point $(j_3+s_k(t_3-t),t)\in {\bar L_k}$ of $P_{k+1}$.
Taking into account the assumption on $dG$ and Lemma~\ref{sum}, we
have

\qtnl{3}
{\bar B_t}= \sum _{i_1+\cdots+i_p=t}
{1 \over i_1!\cdots i_p!} {\partial^t Q_k \over \partial v_1^{i_1} \cdots \partial
v_p^{i_p}} {\bar v_1}^{i_1}\cdots {\bar v_p}^{i_p}
\eqtn
\noindent
Therefore, ${\bar
B_t}=h{\hat B_t}$ where $\hat B_t$ can be treated as a homogeneous
polynomial in ${\bar v_1},\dots , {\bar v_p}$ of the degree $t$
with the coefficients being differential polynomials in $f_1,\dots,f_k$.
Let $t_0$ be the minimal $t$ such that ${\bar B_t}\neq 0$. Then
$t_0\geq 1$ because $Q_k(f_k)=0$, and $t_0\leq t_3$ because ${\bar
B_{t_3}}$ is obtained from $B_{t_3}$ by means of replacing $v_i$
for ${\bar v_i}, 1\leq i \leq p$. One can view $t_0$ as a kind of
multiplicity of the solution $f_k$ in $Q_k$.

\begin{lemma}\label{multiplicity}
$t_0\leq t_4+ {(t_3-t_4)q_{k-1}\over q_k}$.

\end{lemma}

\proof Suppose the contrary. First we observe that the gap between the ordinates
of any pair of consecutive points on $L_k$ is at least $q_k/q_{k-1}$ and that
$e=(t_3-t_4)q_{k-1}/q_k$ is an integer (cf. \cite{Walker}). Hence $L_k$ contains
at most $e+1$ points. Without loss of generality for the sake of conveniency of
notations we assume that $L_k$ contains exactly $e+1$ points (some among them, perhaps,
with zero coefficients $B_t$).

Due to the supposition and the choice of $t_0$ we have ${\bar B_t}=0$ for
$t_4\leq t \leq t_4+e$, i.~e. all the derivatives
$${\partial^t Q_k \over \partial v_1^{i_1} \cdots \partial
v_p^{i_p}}$$
\noindent
of the order $t$ vanish. Fix for the time being non-negative integers $j_1,\dots,j_p$
with the sum $j_1+\cdots+j_p=t_4$. Then
$$0=\sum_{i_1\geq j_1,\dots, i_p\geq j_p; i_1+\dots+i_p=t} {(t-t_4)!\over
(i_1-j_1)!\cdots (i_p-j_p)!} {\partial^t Q_k \over \partial v_1^{i_1} \cdots \partial
v_p^{i_p}} v_1^{i_1-j_1}\cdots v_p^{i_p-j_p}$$
$$=\sum_{l\geq t} {(l-t_4)!\over (l-t)!} {\partial^{t_4} B_l \over \partial v_1^{j_1} \cdots \partial
v_p^{j_p}}$$
\noindent
due to the Euler's formula. The latter equalities can be treated as a linear $(e+1)\times (e+1)$
system with a non-singular matrix. Its non-singularity is justified by the following result
\cite{Lakshman}: if $n_1>\cdots > n_r\geq 0; m_1>\cdots > m_r\geq 0; n_1\geq m_1,\dots, n_r\geq m_r$
then the $r\times r$ matrix with the entries ${n_i \choose m_j}$ is non-singular. Therefore,
$${\partial^{t_4} B_l \over \partial v_1^{j_1} \cdots \partial
v_p^{j_p}} =0$$
\noindent
for any $l$ and any $j_1,\dots, j_p$ with $j_1+\cdots+j_p=t_4$, in particular $B_{t_4}$ vanishes identically, the obtained contradiction
proves the lemma.
\bull

\begin{corollary}\label{change}
If $t_0=t_3$ then the denominator $q_{k-1}=q_k$ does not change.
\end{corollary}

Now we are in position to continue the recursive step of the
procedure constructing $s_{k+1}, f_{k+1}$.
The polygon $P_{k+1}$ either contains the edge with the slope
$s_k$ and with the ordinates $t_0<t_3$, respectively, of its
endpoints, or the edge of $P_{k+1}$ with its above endpoint
$(j_3,t_3)$ has the slope less than $s_k$. In the first case as a
leading edge $L_{k+1}$ one takes an edge of $P_{k+1}$ having a
positive slope $s_{k+1}$ with the ordinate $t_5$ of its upper
endpoint $(j_5,t_5)$ less or equal to $t_0$. In this case
$(j_5,t_5)$ plays the role of a new pivot with $t_5$ being the
multiplicity of $L_{k+1}$. As above one produces the
leading differential polynomial $Q_{k+1}(f_{k+1})$ and as
$f_{k+1}$ chooses a solution of the equation $Q_{k+1}(f_{k+1})=0$.
In the second case the
denominator $q_k=q_{k-1}$ does not increase due to
Corollary~\ref{change}, and as $L_{k+1}$ one takes an edge of
$P_{k+1}$ having a positive slope $s_{k+1}$ with the ordinate
$t_5$ of its upper endpoint (the pivot) $(j_5,t_5)$ less or equal
to $t_3$. The rest is similar to the first case.

Thus, we have described a recursive procedure constructing
$1>s_2>s_3>\cdots$ and $f_1,f_2,f_3,\dots$ which one can view as a
tree.

\begin{lemma}\label{denominator}
i) The common denominator $q$ of $s_2,s_3,\dots$ does not exceed
$2^{m-1}$;

ii) there exists a branch of the tree in which the common denominator $q$
is less or equal to $m$;

iii) every branch of the tree after at most of $q$ steps arrives to
a leading edge with a non-positive slope.
\end{lemma}

\proof First we recall that the multiplicity of any leading edge in $P_2$ is less or
equal to $m$. Therefore, i) follows from  Lemma~\ref{multiplicity}: if at a certain
step the common denominator $q_{k-1}$ is multiplied by
$q_k/q_{k-1}$ then the multiplicity decreases at least by
$q_k/q_{k-1}-1$. After the multiplicity reaches 1, the
denominator does not change anymore.

ii) Let us take at each step of the described recursive procedure the leading edge with
the least possible slope, while the latter is positive. The ordinate of the lower endpoint
of this edge $t_4=0$. Therefore, Lemma~\ref{multiplicity} entails that $t_0\leq t_3q_{k-1}/q_k$,
this implies ii).

iii) follows from  Definition~\ref{symbol} because $k_0\leq q$.
\bull

Assume now that $P_{k+1}$ in the described
procedure contains an edge having a non-positive slope (see Lemma~\ref{denominator}
ii)). Take such edge $L=L_{k+1}$ with the largest possible non-positive slope in $P_{k+1}$.
We have shown above that the coefficient ${\bar B_{t_5}}$ at the
pivot $(j_5,t_5)$ of $L_{k+1}$ equals to $h{\hat B}$ where ${\hat
B}$ is a suitable homogeneous polynomial of the degree $t_5$ in
$d_xf_{k+1}, d_yf_{k+1}$ with the coefficients being differential
polynomials in $f_1,\dots,f_k$. Denote by $\bar B$ the coefficient
at the point $(j_5,0)$ of $P_{k+1}$, being a linear homogeneous operator
in $h$ (one can show that the order of $\bar B$ does not exceed
$t_5$ in the same manner as it was shown that ${\bar B_{t_5}}$ has
the order 0 in $h$). If $\hat B$
contains a term $b(d_xf_{k+1})^l(d_yf_{k+1})^{t_5-l}$ for some $l$
and $b$ being a differential polynomial in $f_1,\dots, f_k$,
then $\bar B$ contains the term $b((d_x)^l(d_y)^{t_5-l}h)$ due to
Lemma~\ref{sum}. Hence the order of $\bar B$ is greater or equal
to $t_5$ (actually, equals $t_5$ as we have seen, although we use below
only that the order of $\bar B$ is positive). In particular,
the slope of $L$ equals 0, and $P_{k+1}$ contains no edges with
negative slopes. In the construction under description $f_{k+1}$ does not appear and
as $h\in F$ we take a solution of the linear homogeneous differential
equation ${\bar B}(h)=0$ (which can be viewed as a leading
equation on $h$).

This completes the construction of the first summand $hG$ of the
solution $H$ of the form (\ref{1}). To obtain the next coefficient
$h_1$ of $H$ we observe that in the expansion of
$T(h_1G^{(-1/q)})$ in the fractional derivatives $\{G^{(i/q)}\}_{-\infty < i < \infty}$ the
highest non-zero term equals ${\bar B}(h_1)G^{(j_5-1/q)}$, taking
into account that this expansion is obtained by means of the shift by $-1/q$ of the
expansion of $T(hG)$ while replacing $h$ for $h_1$. Therefore, for
$h_1\in F$ we get a linear partial differential equation
(not necessary, homogeneous) of the form ${\bar
B}(h_1)=\bar f$ (so, of the same order $t_5$) for an appropriate
${\bar f}\in F$ being a differential
polynomial in $h,f_1,\dots,f_{k_0}$ (in the above notations
$k_0=k$). In a similar way one obtains consecutively
$h_2,h_3,\dots$.

Summarizing, the following theorem is proved.

\begin{theorem}\label{theorem}
Any linear partial differential equation $T=0$ of an order $n$ (see (\ref{1.5}))
for each linear factor $(a_1d_xf_1+a_2d_yf_1)$ of a multiplicity
$m$ of its symbol $symb(T)$ has a non-zero fractional-derivatives series solution
of the form (\ref{1}) with the denominator $q\leq m$.

One can continue every branch of the tree of the described procedure constructing
$1>s_2>s_3>\cdots$ and
$f_1,f_2,\dots$ to a solution of the form (\ref{1}) of $T=0$, and every solution
 of the form (\ref{1}) constructed by a described procedure has the
denominator $q\leq 2^{m-1}$.
\end{theorem}

\begin{corollary}
If an LPDO $T$ has no fractional-derivatives series solutions of the form (\ref{1})
corresponding to a factor $a_1d_xf_1+a_2d_yf_1$ of a multiplicity $m$ of the symbol
$symb(T)$ with the denominator $q<m$ then $T$ is irreducible in $F[d_x,d_y]$. In
particular, if $T$ of an order $n$ has no fractional-derivatives series solutions
with the denominator $q<n$ then $T$ is irreducible.
\end{corollary}

\begin{remark}
The bound $q\leq m$ is sharp as shows the following example. Take $T$
(see (\ref{1.5})) such that $(a_1d_xf_1+a_2d_yf_1)$ has the multiplicity $m$
in $\bar T_n$, the multiplicity greater or equal to $m-i$ in $\bar T_{n-i}$
for every $1\leq i\leq m-2$ and the multiplicity 0 in $\bar T_{n-m+1}$,respectively.
Then the polygon $P_2$ has the edge with the endpoints $(n-m,m)$ and $(n-m+1,0)$
which being taken as a leading one (actually, there is no other choice for a
leading edge), provides the slope $s_2=1/m$.
\end{remark}

\begin{remark}\label{algebraic}
Theorem~\ref{theorem} states the bound $q\leq m$ for a particular solution.
It  is unclear how sharp is the bound $q\leq 2^{m-1}$ for all constructed solutions. The natural question is
whether one can improve it by $m$ (one can verify it for $m\leq 7$ by the
direct calculations)? This would be similar to
 the algebraic situation in which such a  bound on the common
denominator in all Puiseux series (\ref{0}) is well known (see e.~g.
\cite{Walker}). We also mention that for solutions in the canonical form basis \cite{Wasow}
 of linear {\it ordinary} differential equations a
similar to the algebraic situation bound on the common denominator (of the rational
exponents) was established in \cite{G90}.
\end{remark}

\section{Multiplicity of generic fractional-derivatives series solutions}
\label{generic}

In the described recursive construction $f_k$ was chosen as a
solution of the equation $Q_k(f_k)=0$. Different choices of $f_k$
could yield different polygons $P_{k+1}$. Therefore, the set of
(even canonical fractional-derivatives series, see
Remark~\ref{canonical}) solutions of the equation $T=0$ is quite
vast. An interesting open question is whether it is possible to
introduce a concept of a multiplicity of a set of
fractional-derivatives series solutions and relate it to $m$? In the
present section we give a partial answer to this question for the so-called
generic solutions.

We view $Q_k$ as a polynomial in two variables $v_1=d_x f_k, v_2=
d_y f_k$. Note that this polynomial is not homogeneous, consider its
factorization $Q_k=\beta_1^{m_1}\cdots \beta_l^{m_l}\beta$ over $F$
where $\beta$ is homogeneous and $\beta_1,\dots, \beta_l$ are
irreducible non-homogeneous. In the recursive construction from
Section~\ref{solution} we distinguish a case which we call {\it
generic}, namely, when $\beta_i(d_xf_k, d_yf_k)=0$ for a certain
$1\leq i\leq l$ such that $m_i=\min\{m_1,\dots , m_l\}$, and the
point $(d_xf_k, d_yf_k)$ is a non-singular one of the plane curve
$Q_k=0$. In the generic case for the multiplicity of $f_k$ we have
$t_0=m_i$ due to (\ref{3}). One can assign the multiplicity $t_0$ to
the set of all $f_k$ satisfying the generic case. We call a solution
(\ref{1}) generic if for each of $f_2,\dots, f_{k_0}$ the generic
case happens in the construction of (\ref{1}). When $k_0=1$ we call
(\ref{1}) generic as well. At the end of developing any generic
solution we arrive to a polygon $P_{k+1}$ having a leading edge
$L_{k+1}$ with the slope 0. Let the upper endpoint (pivot) of
$L_{k+1}$ be $(j_5,t_5)$, then to this generic solution we assign
the multiplicity $t_5$. Observe that we have assigned the
multiplicity to the set of all the generic solutions (\ref{1}) which
follow the same branch in the tree of the construction from
Section~\ref{solution}.

\begin{proposition}
Any linear partial differential equation $T=0$ of an order $n$

i) has a generic solution of the form (\ref{1});

ii) the sum of multiplicities of the generic solutions does not exceed $n$;

iii) the denominator of every generic solution is less than $n^{O(\log n)}$.
\end{proposition}

\proof Each $Q_k, k\geq 2$ is non-homogeneous, that is why i) is justified
taking into account Theorem~\ref{theorem}.

ii) follows (similar to the algebraic Newton-Puiseux series
\cite{Walker}) by inverse induction along the tree of the procedure described in
Section~\ref{solution} due to the inequality $m_1+\cdots + m_l \leq t_3-t_4$.

The latter inequality together with Lemma~\ref{multiplicity} imply that
$t_0\leq t_3 {q_k \over 2q_k-q_{k-1}}$. Therefore, in developing a generic solution
by means of the procedure from Section~\ref{solution} there are at most
$\log _{3/2} n$ steps at which the denominator augments. At each such step the
denominator grows less than in $n$ times (cf. the proof of Lemma~\ref{multiplicity}),
this entails iii). \bull

\begin{remark}\label{one}
In a particular
case $m=1$ we have $q=k_0=1$, all the solutions of the form (\ref{1})
are canonical, the polygon $P_2$ contains a single edge with a
slope less than 1, namely, the edge with the endpoints $(n-1,1)$
and $(n-1,0)$ having the slope 0. It provides the leading
linear
equation on $h$ of the first order, the leading equation
$(a_1d_x+a_2d_y)f_1=0$ on $f_1$ is linear and of the first order as well,
thus, the multiplicity 1 is assigned to the set of (generic)
solutions in case $m=1$.
\end{remark}

\begin{remark}
While in (\ref{1}) we consider series with decreasing orders of derivatives of $G$, one can
easily verify that an equation $T=0$ for an arbitrary $f\in F$ has a solution of the form
$$\sum _{0\leq i< \infty} h_i G^{(i)}$$
\noindent
where $G=G(f)$, with increasing orders of derivatives of $G$. Thus, continuing the analogy with
plane curves, the latter series could be viewed as corresponding to expanding at finite
points all being regular (so, without proper fractional derivatives, i.~e.
$k_0=1$ in Definition~\ref{symbol}), while (\ref{1})
corresponds to expanding at the infinity.
\end{remark}

\section{Fractional-derivatives series solutions of non-holonomic \\
$D$-modules}\label{ideal}

First let $J=\langle p_1,\dots, p_l\rangle \subset F[d_x,d_y]$ be
a differential (non-holonomic) left ideal of the differential type 1 \cite{Kolchin,
Pankratiev}. This means that the Hilbert-Kolchin polynomial
$K_J(z)=ez+e_0$ of $J$ has the degree 1. Denote by $symb(J)
\subset F[d_xf_1,d_yf_1]$ a homogeneous ideal generated by the
symbols of elements of $J$ (cf. Section~\ref{solution}). Then
$K_J$ coincides with the Hilbert polynomial $K_{symb(J)}$
\cite{Bjork, Sturmfels} (one can also deduce this from the Janet
base of $J$ \cite{schwarz, GS05}, we mention that the concept of
Janet bases was a differential historical predecessor of the one of
Groebner bases). Denote
$g=GCD(symb(J)) \in F[d_xf_1,d_yf_1]$.

\begin{lemma}\label{gcd}
The degree $e$ of the ideal $symb(J)$ coincides with $deg(g)$.
\end{lemma}

\proof Since $symb(J)\subset \langle g \rangle$ it suffices to
verify that $dim_F(\langle g \rangle / symb(J))< \infty$.
Nullstellensatz entails that $(symb(J)/g)\supset
(d_xf_1,d_yf_1)^s$ for a suitable $s$, therefore, the homogeneous
component
$$\langle g \rangle _{deg(g)+s} = g \cdot (d_xf_1,d_yf_1)^s
\subset g \cdot (symb(J)/g)=symb(J) \hfill \quad \bull$$

\begin{remark}\label{holonomic}
If ideal $J\subset F[d_x,d_y]$ is holonomic (so, with differential
type $0$) then $GCD(symb(J))=1$.
\end{remark}

The degree $e$ (being the leading coefficient of the Hilbert-Kolchin
polynomial) is called the typical differential dimension of $J$
\cite{Kolchin, Pankratiev}.

For any homogeneous polynomial $g_0\in F[d_xf_1,d_yf_1]$ and $a\in F$
denote by $mult_a (g_0)$ the multiplicity of the linear form
$d_xf_1+ad_yf_1$ in $g_0$. Also for any $p\in F[d_xf_1,d_yf_1]$ we
denote for brevity $mult_a(p)=mult_a(symb(p))$. W.l.o.g. assume that
$d_yf_1$ does not divide $g$ (otherwise, one can perform a suitable
$C$-linear transformation of $d_x,d_y$). We have $mult_a(symb(J))=
mult_a(g)$ and $e=\sum_a mult_a(g)$ (cf. Lemma~\ref{gcd}). For the
time being fix $a\in
F$ such that $mult_a(g)\geq 1$.

Now we introduce the ring $R=F[d_x,d_y](F[d_y])^{-1}$ of {\it
partial-fractional differential operators} \cite{G05}. Its elements
has the form $p_0b^{-1}$ where $p_0\in F[d_x,d_y], \, b\in F[d_y]$.
One can verify (see \cite{G05}) that $R$ is an Ore ring \cite{Bjork},
any element of $R$ can be written in a form ${\bar b}^{-1}{\bar p}$
for appropriate $ {\bar p}\in F[d_x,d_y], {\bar b}\in F[d_y]$. Thereby
$R= (F[d_y])^{-1} F[d_x,d_y]$ and $p_0b^{-1}= {\bar b}^{-1}{\bar p}$
if and only if ${\bar b}p_0={\bar p}b$. Also in \cite{G05} one can
find the algorithms for addition and multiplication of elements in
$R$.
Any element from $R$ can be written in the form $b^{-1}\sum_{0\leq i\leq
w} b_id_x^i$ for suitable $b,b_i\in F[d_y]$ (because a finite family
of elements from $R$ has a common denominator which belongs to
$F[d_y]$, see \cite{G05}).

For the time being fix $G=G_{(s_2,\dots,s_k)} (f_1,\dots,f_k)$
such that $(d_x+ad_y)f_1=0, \, d_yf_1\neq 0$ (cf.
Section~\ref{series}). Denote by $V=V_G$ the $F[d_x,d_y]$-module which
consists of all fractional-derivatives series of the form
(\ref{1}) added by 0.

\begin{lemma}\label{division}
$V$ is an $R$-module
\end{lemma}

\proof For any $0\neq H\in V$ and $0\neq b\in F[d_y]$ we claim
that $bH\neq 0$. Indeed, let
$$H=hG^{(s)}+\sum_{i\geq 1}
h_iG^{(s-i/q)}, h\neq 0; b=t_nd_y^n+\sum_{0\leq i\leq n-1}
t_id_y^i, t_n\neq 0,$$
\noindent
then
$$bH=ht_n(d_yf_1)^nG^{[s+n)}+\sum_{i\geq 1} {\hat
h_i}G^{(s+n-i/q)} \neq 0$$

For any $H_1\in V$ we need to prove the existence of ${\bar H}\in
V$ such that $b^{-1}H_1={\bar H}$, i.~e. $H_1=b{\bar H}$ (the claim above
implies that $\bar H$ is unique). Let
$H_1=h_{1,0}G^{(s)}+h_{1,1}G^{(s-1/q)}+\cdots; h_{1,0}\neq 0$. Then we
look for ${\bar H}={\bar h}G^{(s-n)}+{\bar
h_1}G^{(s-n-1/q)}+\cdots$. Comparing the coefficients of $H_1$ and
$b{\bar H}$ at $G^{(s)}$, we get $h_{1,0}={\bar h}t_n(d_yf_1)^n$
which yields $\bar h$. Comparing the coefficients at $G^{(s-1/q)}$
yields $\bar h_1$ and so on.
\bull

\begin{remark}\label{isomorphism}
By the same token multiplying by $0\neq p\in F[d_x,d_y]$ on $V$ is
an isomorphism, provided that $(d_x+ad_y)f_1$ does not divide
$symb(p)$.
\end{remark}

The ring $R$ is left-euclidean (as well as right-euclidean) with
respect to $d_x$ over the skew-field $F[d_y](F[d_y])^{-1}$, cf.
Lemma 1.3 \cite{G05}. Hence the ideal $\overline {J}=\langle
p_1,\dots, p_l\rangle \subset R$ is principal, let $\overline
{J}=\langle p \rangle$ for an appropriate $p\in J \subset
F[d_x,d_y]$. Then for any $p_0\in J$ (actually, moreover for $p_0\in
\overline J$) the equalities

\begin{equation}\label{4}
\overline {p_0}p=\overline {b_0}p_0, \overline {b}p=\sum_{1\leq
j\leq l} \overline {p_j}p_j
\end{equation}
\noindent hold for suitable $\overline {p_j}
\in F[d_x,d_y]; \, 0\neq \overline {b},\overline {b_0}\in F[d_y]$.

According to (\ref{4}) we have $symb(\overline
{p_0})symb(p)=symb(\overline{b_0})symb(p_0)$, whence
$$mult_a(p_0)=mult_a(\overline{b_0})+mult_a(p_0)=mult_a(\overline
{p_0})+mult_a(p)\geq mult_a(p),$$ \noindent therefore,
$mult_a(p)\leq mult_a(J)$ since
$$mult_a(J)=mult_a(g)=\min_{p_0\in J} mult_a(p_0).$$

On the other hand, from (\ref{4}) we get
$$mult_a(p)=mult_a(\overline {b})+mult_a(p)=mult_a(\sum_{1\leq j\leq l}
\overline {p_j}p_j)\geq mult_a(J)$$ \noindent Thus, the following
lemma is proved.

\begin{lemma}\label{local}
For any $a\in F$ we have $mult_a(p)=mult_a(J)$.
\end{lemma}

\begin{proposition}\label{localization}
A fractional-derivatives series $H$ with $d_yf_1\neq 0$ (see
(\ref{1})) is a solution of the linear partial differential equation
$p=0$
if and only if $H$ is a solution of
the ideal $J$.
\end{proposition}

\proof If $pH=0$ then from (\ref{4}) we have $0=\overline
{p_0}pH=\overline {b_0}p_0H$. Hence $p_0H=0$ due to
Lemma~\ref{division}. The inverse statement follows again from
(\ref{4}) and Lemma~\ref{division}. \bull

\begin{corollary}\label{solution_ideal}
For any $a_1,a_2\in F$ such that $mult_{(a_1d_xf_1+a_2d_yf_1)}
(symb(J)) \geq 1$, the ideal $J\subset F[d_x,d_y]$ has a solution
of the form (\ref{1}) with a denominator
$q\leq mult_{(a_1d_xf_1+a_2d_yf_1)}(symb(J))$
and
$a_1d_xf_1+a_2d_yf_1=0, grad(f_1)\neq 0$.
\end{corollary}

\proof It follows from Lemma~\ref{local},
Proposition~\ref{localization}, Lemma~\ref{denominator}
and Theorem~\ref{theorem}.
\bull

\begin{remark}
If for every $a\in F$ the ideal $J$ has the multiplicity
$mult_a(J)\leq 1$ then all the solutions of $J$ of the form
(\ref{1}) are canonical, and $J$ has precisely $e=\sum_{a\in
F}mult_a(J)$ (which equals the typical differential dimension of
$J$, cf. Lemma~\ref{gcd}) families of fractional-derivatives
series solutions. Moreover, to each of these families a
multiplicity 1 can be naturally assigned (cf. Remark~\ref{one}).
\end{remark}

Finally, let $U\subset (F[d_x,d_y])^l$ be a (non-holonomic) $F[d_x,d_y]$-module of
the differential type at least 1, obviously, the differential type
does not exceed 2 (recall that the differential type equals the
degree of the Hilbert-Kolchin polynomial of $U$ \cite{Kolchin,
Pankratiev}). Denote by $u_1,\dots,u_l$ a free base of
$(F[d_x,d_y])^l$. For any $1\leq r\leq l$ consider the submodule
$U_r=\{\sum_{r\leq i\leq l} p_iu_i \in U\}$ where $p_i\in
F[d_x,d_y]$. Denote by $J_r=\{p_r\}\subset F[d_x,d_y]$ the left
ideal being the projection of $U_r$ on the $r$-th component. Then
the differential type of $U$ coincides with the maximum of the
differential types of $\{J_r\}_{1\leq r\leq l}$ (one can verify
this, e.~g. using the Janet bases of $\{J_r\}_{1\leq r\leq l}$
which provide a triangular Janet base of $U$). Take the minimal
$r_0$ such that $J_{r_0}$ has the differential type at least 1.

One has the natural action $U\times V^l \rightarrow V$ on the free
$F[d_x,d_y]$-module $V^l=\{\sum_{1\leq i\leq l} H_iv_i\}$ where
$H_i\in V$ (cf. Lemma~\ref{division}) and $v_1,\dots, v_l$ is a
free base of $V^l$. If $\sum_{1\leq i\leq l}p_iH_i=0$ then we call
$\sum_{1\leq i\leq l}H_iv_i$ a solution of $\sum_{1\leq i\leq
l}p_iu_i$ (we shall choose $G$ and thereby, $V=V_G$ later). We are
looking for a solution of the form $\sum_{1\leq i\leq l}H_iv_i$ of
the module $U$.

First we put $H_{r_0+1}=\cdots=H_l=0$ and as $H_{r_0}\neq 0$ take
a fractional-derivatives series being a solution of the ideal
$J_{r_0}$ according to Corollary~\ref{solution_ideal} in case when
the differential type of $J_{r_0}$ equals 1. When the differential
type of $J_{r_0}$ equals 2, in other words, $J_{r_0}=0$, we take
as $H_{r_0}\neq 0$ an arbitrary fractional-derivatives series. In
both cases $H_{r_0}v_{r_0}$ is a solution of the submodule
$U_{r_0}$. Thus, we have chosen $G=G_{(s_2,\dots,s_k)}
(f_1,\dots,f_k)$ and thereby, $V=V_G$. As above we can assume w.l.o.g.
that in the equation $(a_1d_x+a_2d_y)f_1=0$ we have $a_1\neq 0$,
so $(d_x+ad_y)f_1=0$ (performing if necessary a suitable
$C$-linear transformation of $d_x,d_y$).

Now we construct $H_r$ by recursion on $r_0-r\geq 0$. Suppose that
we have already constructed an element $\sum_{r+1\leq i\leq l}
H_iv_i$ being a solution of $U_{r+1}$ for some $r+1\leq r_0$.
Since $J_r$ has the differential type 0 (due to the choice of
$r_0$), $J_r$ contains a certain element $0\neq b\in F[d_y]$.
Consider a corresponding element $u=bu_r+\sum_{r+1\leq i\leq l}
p_iu_i \in U_r$. According to Lemma~\ref{division} one can find
$H_r\in V$ such that $bH_r+\sum_{r+1\leq i\leq l} p_iH_i=0$. For
any element ${\bar u}= \sum_{r\leq i\leq l}{\bar p_i}u_i\in U_r$
applying the left euclidean division in $R$ one can represent
${\bar u}={\bar p_r}b^{-1}u+{\hat u}$ for an appropriate ${\hat
u}\in U_{r+1}$. Then ${\hat u}(\sum_{r\leq i\leq l} H_iv_i)=0$ by
the recursive hypothesis. Besides,
${\bar p_r}b^{-1}u(\sum_{r\leq i\leq l} H_iv_i)=0$ because of
Lemma~\ref{division}. Hence
${\bar u}(\sum_{r\leq i\leq l} H_iv_i)=0$ which completes the
recursive step.

Summarizing, the following main theorem of the paper is proved.

\begin{theorem}\label{module}
Any (non-holonomic) module in $(F[d_x,d_y])^l$ of the differential type at least 1
has a fractional-derivatives series non-zero solution.
\end{theorem}

\begin{remark}\label{three}
One could consider an ideal $J\subset F[d_{x_1},\dots,d_{x_t}]$
still of the differential type 1 with a number of
derivatives $t\geq 3$ and ask whether $J$ has always a
fractional-derivatives series solution? The answer to this
question is negative already for $t=3$ and an ideal
$J=\langle p_1,p_2 \rangle \subset F[d_{x_1},d_{x_2},d_{x_3}]$
(being generic of the differential type 1) generated by an
operator $p_1$ of the first order and $p_2$ of the second order.
\end{remark}

\section{Duality between non-holonomic ideals and
fractional-derivatives series solutions}\label{duality}

There is a well-known duality \cite{Kolchin} between (left)
differential ideals and their spaces of solutions (being an analogue
of the duality between radical ideals and varieties in algebraic
geometry). To establish a similar duality for non-holonomic ideals
in $F[d_x,d_y]$ (so, of the differential type $1$) we need to make
use of the equivalence relation on ideals introduced in \cite{GS05}.
We say that non-holonomic ideals $0\neq J,J_0\subset F[d_x,d_y]$ are
{\it equivalent} if the leading coefficients of degree 1 (see
Section~\ref{ideal}) Hilbert-Kolchin polynomials of three ideals
$J,J_0,J\cap J_0$ coincide (denote these leading coefficients by
$e$), then moreover, $GCD(symb(J))=GCD(symb(J_0))=GCD(symb(J\cap
J_0))$ and the degree of the latter polynomial equals $e$ (see
Lemma~\ref{gcd}). In this case ideal $\langle J,J_0\rangle$ is also
non-holonomic and $GCD(symb(\langle J,J_0\rangle))=GCD(symb(J))$ as
well \cite{GS05, Cassidy, Sit} and moreover, clearly four ideals $J,
\, J_0, \, J\cap J_0, \, \langle J,J_0\rangle$ are equivalent.
Equivalence classes of ideals play a similar role to classes (in
algebraic geometry) of plane curves with the same sets of
$1$-dimensional components.
In this Section we
prove that the sets of fractional-derivatives series solutions of
equivalent non-holonomic ideals coincide and that there is a duality
between the equivalence classes of ideals and their respective sets
(which basically means that to distinct classes correspond distinct
sets).

In this Section we keep the notations from Section~\ref{ideal}. The
next lemma states that the multiplication by $0\neq p\in F[d_x,d_y]$
on $D$-module $V=V_G$ is an epimorphism (the conditions on its
injectivity follow from Theorem~\ref{theorem}, cf. also
Remark~\ref{isomorphism}).

\begin{lemma}\label{epi}
For any $0\neq p\in F[d_x,d_y]$ we have $pV=V$.
\end{lemma}

{\bf Proof}. Let $0\neq H_1=\overline{h}G^{(s_0)}+\sum_{i\ge 1}
\overline{h_i}G^{(s_0-i/q)} \in V$ (see (\ref{1})). We search for
$H=hG^{(s)}+\sum_{i\ge 1} h_iG^{(s-i/q)} \in V$  such that $pH=H_1$.
Treating $h,s$ as indeterminates we get from Lemma~\ref{sum} that

\qtnl{51} p(hG^{(s)})=\eta G^{(s+\kappa /q)}+\sum_{1\le i\le \kappa}
\eta _i G^{(s+\kappa /q -i/q)} \eqtn

\noindent for certain integer $\kappa \ge 0$ and linear ordinary
differential operators $\eta \neq 0, \eta _i$ in $h$ (with
coefficients being differential polynomials in $f_1,\dots , f_k$
which we recall are assumed to be fixed). Indeed, such $\kappa$ with
non-zero $\eta$ exists since in the expansion of $p(hG^{(s)})$ the
coefficient at $G^{(s)}$ in (\ref{51}) equals $p(h)$, so is a
non-zero linear ordinary differential operator. Therefore, we put
$s=s_0-\kappa /q$, and there exists $h\in F$ for which
$\eta(h)=\overline h$ (because $F$ is differentially closed). At the
next step comparing the coefficients of $pH$ and $H_1$ at
$G^{(s+\kappa /q-1/q)}$ one can find $h_1$ from an equation of the
form $\eta(h_1)=\tilde h$ for certain $\tilde{h} \in F$, and so on
one can find $h_2, h_3, \dots$ consecutively. \bull \vspace{4mm}

First consider two equivalent non-holonomic ideals $J,J_0\subset
F[d_x,d_y]$. We claim that four sets of all fractional-derivatives
series solutions of $J,J_0, J\cap J_0, \langle J,J_0\rangle$,
respectively, coincide. As in Section~\ref{ideal} one can suppose
w.l.o.g. that $d_y f_1$ does not divide $g=GCD(symb(J))$ and
consider left ideals $\overline{J}, \overline{J_0} \subset R$ being
principal. Let $\overline{J}=\langle p\rangle, \, \overline{J\cap
J_0} = \langle p_0 \rangle$ for suitable generators $p,p_0\in
F[d_x,d_y]$. Then $bp_0=p_2p$ for appropriate $p_2\in F[d_x,d_y], \,
0\neq b\in F[d_y]$. Lemma~\ref{local} implies that $symb(p_0)$
coincides with $symb(p)$ up to a power of $d_y f_1$, whence
$symb(p_2)$ is a power of $d_yf_1$. Therefore,
Proposition~\ref{localization} and Remark~\ref{isomorphism} entail
the required claim on coincidence of the sets of
fractional-derivatives series solutions of ideals $J, \, J\cap J_0$
(and in a similar way also of $J_0, \, \langle J,J_0\rangle$).

Now let non-holonomic ideals $J\subsetneq J_1 \subset F[d_x,d_y]$ be
non-equivalent. Our purpose is to find a fractional-derivatives
series solution of $J$ being not a solution of $J_1$. Denote
$g_1=GCD(symb(J_1))$, clearly $g_1|g$. Again we suppose w.l.o.g.
that $d_y f_1$ does not divide $g$. Let $\overline{J}=\langle p
\rangle, \, \overline{J_1}=\langle p_1 \rangle \subset R$. Then
$b_1p=p_3p_1$ for suitable $p_3\in F[d_x,d_y], \, 0\neq b_1\in
F[d_y]$. Since $J,J_1$ are not equivalent we have $\deg (g_1)<\deg
(g)$ (see Lemma~\ref{gcd}) and because of that Lemma~\ref{local}
implies that $symb(p_3)$ has a divisor of the form
$d_xf_1+a_3d_yf_1$ for a certain $a_3\in F$. Due to
Theorem~\ref{theorem} there exists a fractional-derivatives series
solution $H_1=\sum_{i\ge i_0} h_iG^{(-i/q)}$ (see (\ref{1})) of
equation $p_3H=0$ for appropriate $G=G_{s_2,\dots,s_k}
(f_1,f_2,\dots, f_k)$ where $d_xf_1+a_3d_yf_1=0$. Now we apply
Lemma~\ref{epi} to $p_1$ and obtain a fractional-derivatives series
$H\in V_G$ such that $p_1H=H_1$. Therefore, $b_1pH=0$ and hence
$pH=0$ in view of Lemma~\ref{division}. Thus, $H$ is a desired
solution of $J$ being not a solution of $J_1$.

Finally, consider non-equivalent non-holonomic ideals $J,J_1\subset
F[d_x,d_y]$ and assume that their respective sets of
fractional-derivatives series solutions coincide. Then ideal
$\langle J,J_1\rangle$ has also the same set of
fractional-derivatives series solutions, in particular $\langle
J,J_1\rangle$ is non-holonomic by virtue of Remark~\ref{holonomic}
and of Theorem~\ref{theorem}. Therefore, due to the proved above
three ideals $J, \, \langle J,J_1\rangle, \, J_1$ are equivalent
which contradicts to the assumption.

We summarize the proved duality in the following

\begin{proposition}\label{dual}
Non-holonomic left ideals $J,J_0 \subset F[d_x,d_y]$ are equivalent
if and only if they have the same sets of all fractional-derivatives
series solutions.
\end{proposition}

For a non-holonomic ideal $I\subset F[d_x,d_y]$ denote by $[I]$ the
equivalence class of non-holonomic ideals which contains $I$ and by
$V(I)$ the set of all fractional-derivatives series solutions of
$I$. In \cite{GS05} we define the following partial ordering on the
classes: $[J]$ is {\it subordinated} to $[I]$ if there exist ideals
$J_1\in [J], \, I_1\in [I]$ such that $J_1\subset I_1$.

\begin{corollary}\label{subordinated}
$[J]$ is subordinated to $[I]$ if and only if $V(I)\subset V(J)$.
\end{corollary}

{\bf Proof}. Let $V(I)\subset V(J)$, then $V(\langle
I,J\rangle)=V(I)$. Proposition~\ref{dual} entails that $\langle
I,J\rangle$ is equivalent to $I$, hence $[J]$ is subordinated to
$[I]$.

The inverse implication is evident. \bull \vspace{4mm}

Now we connect the subordination relation with localizations of
ideals in the ring $R=F[d_x,d_y](F[d_y])^{-1}$ (see
Section~\ref{ideal}).

\begin{proposition}
If $V(I)\subset V(J)$ then $\overline{J}\subset \overline{I}\subset
R$ (provided that $d_yf_1$ does not divide $GCD(symb(I))$).
\end{proposition}

{\bf Proof}. Let $\overline{J}=\langle p\rangle, \,
\overline{I}=\langle q\rangle$ for suitable $p,q\in F[d_x,d_y]$ (cf.
above).  Then $p=p_0q$ for appropriate $p_0\in R$. Whence
$V(I)\subset V(J)$ relying on Proposition~\ref{localization}. \bull


\vspace{4mm}







For a pair of left ideals $J\subset J_1\subset
F[d_{x_1},\dots,d_{x_m}]$ we have introduced in \cite{GS05} a
concept of relative syzygies. Namely, let $J_1=\langle
p_1,\dots,p_t\rangle$, then we define the left {\it module of
relative syzygies}
$$Syz(J,J_1)=\{(q_1,\dots,q_t): \, \sum_{1\le i\le t} q_ip_i \, \in
J; \, q_i\in F[d_{x_1},\dots,d_{x_m}], \, 1\le i\le t\}.$$ \noindent
Making use of \cite{Quadrat} one can verify \cite{GS05} that module
$Syz(J,J_1)$ is independent of a choice of generators
$p_1,\dots,p_t$. Let us denote by $U(J)\subset F$ the space of
solutions of $J$ which can be treated as a $C$-vector space. It was
proved in \cite{GS05} that the quotient $U(J)/U(J_1)$ is isomorphic
to $U(Syz(J,J_1))\subset F^t$.

Here we establish a similar result for non-holonomic ideals
$J\subset J_1\subset F[d_x,d_y]$ and their spaces of
fractional-derivatives series solutions $V_G(J)\subset V_G$ of the
form (\ref{1}) for any $G$ fixed for the time being (see
Section~\ref{ideal}), again we treat $V_G(J)$ as a $C$-vector space.
As in Section~\ref{ideal} one can assume w.l.o.g. that $d_yf_1\neq
0$.

Mapping $\psi:v\to (p_1,\dots,p_t)^T v$ assures a monomorphism
$V_G(J)/V_G(J_1)\hookrightarrow V_G(Syz(J,J_1))$. To show that it is
an epimorphism take an arbitrary vector $(w_1,\dots,w_t)\in
V_G(Syz(J,J_1))\subset V_G^t$. The following property holds: for any
$q_1,\dots,q_t\in F[d_x,d_y]$ such that $\sum_{1\le i\le t}
q_ip_i=0$ (moreover, one can suppose that $\sum_{1\le i\le t} q_ip_i
\in J$) we have $\sum_{1\le i\le t} q_iw_i=0$. Clearly, this
property holds also for any $q_1\dots,q_t\in R$ (see
Section~\ref{ideal}). Consider principal ideal
$\overline{J_1}=\langle p_1,\dots,p_t\rangle=\langle p\rangle
\subset R$ for suitable $p\in F[d_x,d_y]$. Then $p=b^{-1}\sum_{1\le
i\le t} \eta_ip_i$ for appropriate $0\neq b\in F[d_y], \,
\eta_1,\dots, \eta_t \in F[d_x,d_y]$. Denote $w=\sum_{1\le i\le t}
\eta_iw_i \in V_G$. Due to Lemma~\ref{epi} there exists $v\in V_G$
such that $(\sum_{1\le i\le t} \eta_ip_i)v=w$. For each $1\le i_0\le
t$ one can find $\lambda_{i_0}\in R$ for which
$p_{i_0}=\lambda_{i_0}p$, the mentioned above property implies that
$w_{i_0}=\lambda_{i_0}b^{-1}(\sum_{1\le i\le t} \eta_iw_i)$, hence
$w_{i_0}=p_{i_0}v$, i.~e. $\psi(v)=(w_1,\dots,w_t)$. Finally, we
check that $v\in V_G(J)$. Indeed, an arbitrary $q\in J$ can be
represented as $q=\sum_{1\le i\le t} q_ip_i$ for certain
$q_1,\dots,q_t\in F[d_x,d_y]$, then $(q_1,\dots,q_t)\in Syz(J,J_1)$,
whence $qv=\sum_{1\le i\le t} q_iw_i=0$. Thus, in the introduced
notations we have proved the following

\begin{proposition}
For any non-holonomic ideals $J\subset J_1\subset F[d_x,d_y]$ and
$G$ there is an isomorphism of $C$-vector spaces $V_G(J)/V_G(J_1)$
and $V_G(Syz(J,J_1))$.
\end{proposition}

One can deduce Proposition~\ref{dual} from the latter Proposition
invoking Theorem~\ref{module}.

It would be interesting to clarify, whether for non-holonomic ideals
$I,J\subset F[d_x,d_y]$ the equality
$$GCD(symb(I))\cdot GCD(symb(J))=GCD(symb(\langle I,J\rangle))\cdot
GCD(symb(I\cap J))$$ \noindent holds? Observe that the degrees of
the polynomials in both sides of the latter equality coincide in
view of \cite{Cassidy, Sit} taking into account Lemma~\ref{gcd}. A
more subtle question is whether for any $G$ the equality $V_G(I\cap
J)=V_G(I)+V_G(J)$ is true?

\section{Completeness of fractional-derivatives solutions for separable
linear partial differential operators}\label{completeness}

Let $T=T_n+\cdots +T_0\in F[d_x,d_y]$ be a {\it separable} LPDO, i.e. its symbol $symb(T)=
{\overline {T_n}} = \prod_{1\leq i\leq n} (d_xf-a_id_yf)$ is the product of $n$ pairwise
distinct homogeneous linear forms in $d_xf,d_yf$. One can always bring $symb(T)$ to this form monic with respect to
$d_xf$ making, if necessary, a $C$-linear transformation of $d_x,d_y$ in case when $symb(T)$
has a divisor $d_yf$.

For each $1\leq i\leq n$ the equation $T=0$ has a fractional-derivatives series solution of the form (due to Theorem~\ref{theorem})

\begin{equation}\label{2.1}
h_{0,i}G^{(0)}(f_i)+h_{1,i}G^{(-1)}(f_i)+\cdots
\end{equation}

\noindent
where $d_xf_i-a_id_yf_i=0$ and $h=h_{0,i}$ satisfies the first-order LPDE

\begin{equation}\label{2.2}
{{\overline {T_n}}(f_i) \over (d_xf_i-a_id_yf_i)} (d_xh-a_id_yh)+{\overline {T_{n-1}}}(f_i)h=0
\end{equation}

\noindent
We observe that $h_{j,i};j=1,2,\dots$ satisfy similar to (\ref{2.2}) equations with the
highest (first-order) form ${{\overline {T_n}}(f_i) \over (d_xf_i-a_id_yf_i)} (d_xh_{j,i}-a_id_yh_{j,i})$,
being not necessary homogeneous.

From now on throughout this section we assume that $F$ is the field of meromorphic functions in a certain domain
$M\subset \CC^2$, thus the coefficients of $T$ belong to $F$. For a suitable point $(x_0,y_0)
\in M$ the series (\ref{2.1}) can be rewritten as a formal power series in $x-x_0,y-y_0$. Our
goal is to find a point $(x_0,y_0)$ and look for solutions of $T=0$ as power series in $x-x_0,y-y_0$.

We choose a point $(x_0,y_0)
\in M$ such that all the coefficients of $T$ at this point are defined and in addition, the
values $a_i(x_0,y_0)$ are pairwise distinct for $1\leq i\leq n$. The latter is equivalent to
that the discriminant of  $symb(T)$ does not vanish at this point. Therefore, all the points
of $M$ out of an appropriate analytic subvariety of $M$ of the dimension 1 satisfy these requirements.

One takes a solution $f_i$ (being a power series in $x-x_0,y-y_0$) of the equation
$d_xf_i-a_id_yf_i=0$ with a vanishing free coefficient (which we denote by $f_i(x_0,y_0)=0$) and
with a non-vanishing vector of coefficients
at the first powers of $x-x_0,y-y_0$ (which we denote by $(d_xf_i,d_yf_i)(x_0,y_0)$, thereby $d_yf_i(x_0,y_0)\neq 0$). We observe that this LPDE has always a
 solution with arbitrary chosen free coefficient and non-vanishing vector of the coefficients
at the first powers of $x-x_0,y-y_0$ since the vector of the coefficients $(1,-a_i)$
at its highest (first) derivatives does not vanish at the point  $(x_0,y_0)$. Hence the free
coefficient of the power series $d_xf_i-a_jd_yf_i$ does not vanish when $j\neq i$ due to the
requirement on the discriminant.
Therefore, by the same token one can find a solution $h$ of the equation (\ref{2.2}) with a
non-zero free coefficient which we denote by $h(x_0,y_0)\neq 0$.

We take an arbitrary solution of $T=0$, being a power series in $x-x_0,y-y_0$ and intend to represent it as a sum of $n$ solutions of the form (\ref{2.1}) (for
$1\leq i \leq n$) in
which $G^{(0)}(f_i)$ is replaced by its specialization (see Remark~\ref{specialisation})
$$\sum_{j\geq 0} c_{j,i} {f_i^j \over j!}$$
\noindent
with indeterminate coefficients $c_{j,i}\in \CC$. Then
$$G^{(-l)}=\sum_{j\geq 0} c_{j,i} {f_i^{j+l} \over (j+l)!}.$$

Suppose that by recursion on $k$ the coefficients $c_{j,i}$ for $j\leq k-1, 1\leq i\leq n$ are
already produced. Our purpose is to produce $c_{k,i},1\leq i\leq n$. Clearly, any solution of
$T=0$ being a power series of the form $\sum_{p,q\geq 0} b_{p,q}(x-x_0)^p(y-y_0)^q$ is
determined by the coefficients $b_{p,q}$ with $0\leq p\leq n-1$.

For each $0\leq p\leq n-1$ the contribution of the term at $c_{k,i}$ (see (\ref{2.1})) into
$b_{p,k-p}$ equals to
$$h(x_0,y_0)(d_x^pd_y^{k-p}{f_i^k \over k!})(x_0,y_0)=h(x_0,y_0)(a_i^p(d_yf_i)^k)(x_0,y_0)$$
\noindent
taking into account that $f_i(x_0,y_0)=0, d_yf_i(x_0,y_0)\neq 0, h(x_0,y_0)\neq 0$.

Therefore, we obtain a linear (algebraic) system (in general, not necessary homogeneous) on
$c_{k,i}; 0\leq i\leq n-1$ with the matrix being of the van-der-Monde type $(a_i^p(x_0,y_0))$. This
allows one to find uniquely $c_{k,i}; 0\leq i\leq n-1$ and thereby, carry out the recursive
step.

\begin{theorem}
For a separable LPDO $T$ of the order $n$ with the coefficients being meromorphic in a certain complex domain
$M$ the sum of $n$ spaces of specialisations of
fractional-derivatives series solutions of $T=0$ of the form (\ref{2.1})
(for fixed $f_i$ and $h_{j,i}$) coincides with the space of all the solutions of $T=0$ as
formal power series in $x-x_0,y-y_0$ for any point $(x_0,y_0)$ from $M$ out of a suitable analytic
subvariety of the dimension 1.
\end{theorem}

It would be interesting to extend this theorem to a non-separable LPDO. Let us also mention that in
\cite{GS04} an algorithm for factoring a separable LPDO was produced.

\section{Applications to studying first-order factors of a linear partial differential operator}
\label{applications}

\subsection{Finding first-order factors of a linear partial differential operator}\label{factors}

Let $T=T_n+\dots+T_0$ be a LPDO of an order $n$ in 2 independent variables, where
$T_j=\sum_i a_{i,j-i}d_x^id_y^{j-i}$ is a sum of the derivatives of the order $j$.
We assume that the coefficients
$a_{i,j}$ are taken from the field $\QQ(x,y)$ in order to design algorithms, while $f$ is taken
from a universal field $F$ (cf. Section~\ref{completeness}).

As we are looking for the first-order factors of $T$ of the form $L=d_x+ad_y+b \in F[d_x,d_y]$ we need to study the
solutions of $L=0$ (w.l.o.g. one can assume that the coefficient at $d_x$ of $L$ does not
vanish, otherwise one can change the roles of $x$ and $y$).
Take any solution $f$ of the symbol $(d_x+ ad_y)f=0$ of $L$ such that $d_yf\neq 0$ and consider
$G=G^{(0)}(f)$ (cf. Section~\ref{completeness}).
For any $h\in F$ being a ``particular'' solution of $L=0$,
we have that $hG$ is a fractional-derivatives series solution of $L=0$.

\begin{lemma}\label{correspondence}
An operator $T$ has a right first-order factor $L$ if and only if the equation $T=0$ has a
solution of the form $hG$.
\end{lemma}

{\bf Proof.} If $T$ has a right factor $L$ then $T$ has a solution $hG$.

Conversely, assume that $T=0$ has a solution $hG$. Dividing $T$ with remainder by $L$ one
can represent $T=SL+\sum_{0\leq i\leq n} b_id_y^i$ for a suitable operator $S$. Consider
the largest $k$ such that $b_k\neq 0$. Then in the expansion of $(\sum_{0\leq i\leq k} b_id_y^i)hG$ in $\{G^{(s)}\}$ the coefficient at $G^{(k)}$ equals $b_kh(d_yf)^k\neq 0$. The obtained
contradiction shows that $T=SL$. \bull

Thus, we are looking for a solution of $T=0$ of the form $hG$. Expanding $T(hG)=A_0G^{(0)}+\cdots+A_{n}G^{(n)}$, we get first that $A_n=symb(T)$. Therefore, we fix for the time being a linear divisor
of the form $d_xf+ad_yf$ of $symb(T)$ and assume that this divisor vanishes.
Thereby, the calculations below (arthmetic manipulations and polynomial factoring) will be carried out over the field $\QQ(x,y)[a]$. This can be fulfilled representing $\QQ(x,y)[a]\simeq \QQ(x,y)[z]/(g)$ where $g\in \QQ(x,y)[z]$ is the minimal polynomial of $a$ (see \cite{G86}).
So, we obtain $n$
equations $A_0=\cdots=A_{n-1}=0$ treated as LPDO in $h$ with the coefficients being non-linear
differential polynomials in $f$. We denote the ring of all these polynomials by $P=\QQ(x,y)[a]\{d_xf,d_yf\}$. Applying to $A_0=\cdots=A_{n-1}=0$ the procedure of constructing a Janet base \cite{schwarz} one
gets the conditions of solvability in $h$ of $A_0=\cdots=A_{n-1}=0$ expressed as a disjunction of
systems of the form
\begin{equation}\label{1.1}
p_1=\cdots=p_l=0, p_0\neq 0
\end{equation}
\noindent
where $p_i\in P$. Using the relation $d_xf+ad_yf=0$ one can reduce each $p_i$ to an (ordinary)
differential polynomial $\bar p_i$ in $d_yf$. Denote the ring of ordinary differential polynomials
by $R=\QQ(x,y)[a]\{d_yf\}$.

Applying to the formula ${\bar p_1}=\cdots={\bar p_l}=0, {\bar p_0}\neq 0$ the subroutine of the elimination procedure
in the theory of ordinary differentially closed fields from \cite{S56} (see also \cite{G89} where its improvement with a better complexity bound was designed) one obtains an equivalent disjunction of systems
of the form
\begin{equation}\label{1.2}
r=0, r_0\neq 0
\end{equation}
\noindent
for suitable differential polynomials $r,r_0\in R$. Briefly, this subroutine consists in alternative
executing 2 types of steps while there are more than one equality of (ordinary) differential polynomials.
The first type of steps is executed when all the highest derivatives occurring in these polynomials are
equal, in this case the algorithm calculates their GCD viewing them as (algebraic) polynomials in this
highest derivative (and branching depending on vanishing the leading coefficients). Else, if not all the
highest derivatives are equal, as the second type of steps one can diminish the highest derivative.
Moreover, if $r$ contains the $d_y^kf$ as its
highest derivative then $r$ considered as an (algebraic) polynomial in the ring $K=\QQ(x,y)[a][f,d_yf,\dots,d_y^kf]$ is
irreducible. In addition, $r_0$ is less than $r$ with respect to the term ordering, i.~e. if $r_0$
contains $d_y^{k_0}f$ as its highest derivative then either $k_0<k$ or $k_0=k$ and the degree of
$r_0$ with respect to $d_y^kf$ is less than the similar degree of $r$.

Replace $d_xf$ by $-ad_yf$ in $d_xr$. This yields a differential polynomial ${\hat r} \in R$ of the
order at most $k+1$ (its role is similar to an $S$-pair in Janet type algorithm \cite{schwarz}). If
$\hat r$ does not belong to the differential ideal $\langle r \rangle \subset R$, we again apply to
the system $r={\hat r}=0, r_0\neq 0$ the used above subroutine from the elimination procedure and get
an equivalent disjunction of systems of the form (\ref{1.2}) with less term ordering than of $r$ and
continue as above.

Now assume that $\hat r$ belongs to $\langle r \rangle$. Then we claim that any solution of (\ref{1.2})
provides a solution of (\ref{1.1}). Indeed, otherwise, the ideal $\langle r, d_xf+ad_yf\rangle \subset
P$ would contain an appropriate power $r_0^s$ \cite{Kolchin}, p.146-148. This yields a relation of
the form
$$r_0^s=\sum_{i,j}A_{i,j}d_x^id_y^jr +\sum_{i,j}B_{i.j}d_x^id_y^j(d_xf+ad_yf)$$
\noindent
for suitable $A_{i,j},B_{i,j}\in P$. Replacing in this relation $d_xf$ for $-ad_yf$ and taking into
account that $\hat r$ belongs to $\langle r \rangle$, we deduce that
\begin{equation}\label{1.3}
r_0^s=\sum_j {\hat A_j} d_y^jr
\end{equation}
\noindent
for certain ${\hat A_j}\in R$. From the equation $d_yr=0$ we express
$$d_y^{k+1}f={\hat B_{k+1}}/{\partial r \over \partial (d_y^kf)}$$
\noindent
for an appropriate ${\hat B_{k+1}} \in K$.
After that express
successively
$$d_y^{k+2}f={\hat B_{k+2}}/{\partial r \over \partial (d_y^kf)},
d_y^{k+3}f={\hat B_{k+3}}/{\partial r \over \partial (d_y^kf)},\dots.$$
\noindent
Substitute these expressions in (\ref{1.3}), this results in the  equality
$$r_0^s({\partial r \over \partial (d_y^kf)})^t=Ar$$
\noindent
for some $t$ and $A\in K$.
But $r$ is irreducible in $K$
and $r_0$ is less than $r$ with respect to the term ordering. The obtained contradiction proves the
claim and the following theorem.
\begin{theorem}
There is an algorithm which tests whether an operator $T \in \QQ (x,y)[d_x,d_y]$ has a first-order factor with the
coefficients in a universal field $F$. The algorithm invokes two subroutines: the elimination of an
unknown function in a system of LPDO's (in other words, a parametric Janet base), and a subroutine
from the elimination procedure in the theory of ordinary differentially closed fields.
\end{theorem}

\begin{remark}
If one uses a direct method of finding the coefficients of a first-order operator $L$ and of an
$(n-1)$-th order $Q$ such that $T=QL$, then one has to apply an elimination in the theory of
{\it partial} differentially closed fields whose complexity is unclear how to estimate in a
reasonable way (cf. \cite{S56, G89}).
\end{remark}

\begin{remark}
One can also search for left first-order factors of an LPDO (by means of
considering an adjoint operator).
\end{remark}

\begin{corollary}
There is an algorithm to factor LPDO's of the orders at most 3.
\end{corollary}

\subsection{Intersection of principal first-order ideals}\label{ideals}

In this subsection by $F$ we denote a differential field with derivatives $d_x,d_y$.

First consider the ideals $I_i=\langle d_x+ad_y+b_i\rangle$ with the same highest (first-order)
forms where $a,b_i\in F, 1\leq i\leq n$.

\begin{proposition}\label{intersection}
The ideal $I_1\cap \cdots \cap I_n$ is principal
\end{proposition}

{\bf Proof.} Denote $E=d_x+ad_y$. The ring $F[E]$ is left-euclidean, therefore, the intersection
${\hat {I_1}}\cap \dots \cap {\hat{I_n}} = \langle Q \rangle \subset F[E]$ is principal where we denote
${\hat{I_i}}=\langle d_x+ad_y+b_i \rangle \subset F[E]$ and $Q=q_sE^s+\cdots + q_0$ for certain
$q_0,\dots, q_s\in F, q_s \neq 0$ and $s\leq n$.

Our aim is to prove by induction on $n$ that $I_1\cap \cdots \cap I_n = \langle Q \rangle \subset
F[d_x,d_y]$. Assume that it is already proved and consider the intersection $I_1\cap \cdots \cap
I_n \cap I_{n+1}$. There can occur two cases. Either ${\hat {I_1}}\cap \dots \cap {\hat{I_n}} \cap
{\hat I_{n+1}} = {\hat {I_1}}\cap \dots \cap {\hat{I_n}}$, in this case $Q=N(d_x+ad_y+b_{n+1})$ for a
suitable $N\in F[E]$, therefore, $I_{n+1} \supset \langle Q \rangle$ and $I_1\cap \cdots \cap I_n \cap
I_{n+1}= \langle Q \rangle$.

Or else ${\hat {I_1}}\cap \dots \cap {\hat I_n} \supsetneq
{\hat {I_1}}\cap \dots \cap {\hat{I_n}} \cap
{\hat I_{n+1}} = \langle M \rangle$ for an appropriate $M=m_{s+1}E^{s+1}+\cdots + m_0 \in F[E]$ with
$m_0,\dots , m_{s+1} \in F$. Clearly, $M\in I_1\cap \cdots \cap
I_n \cap I_{n+1}$. It is necessary to show that for any $V\in I_1\cap \cdots \cap
I_n \cap I_{n+1}$ we have $V\in \langle M \rangle$. Since the highest derivative with respect to
$d_x$ which occurs in $M$ is $d_x^{s+1}$, one can divide $V$ by $M$ with remainder and get
$V=WM+U$ where  $W,U\in F[d_x,d_y]$ for a certain $U\in I_1\cap \cdots \cap
I_n \cap I_{n+1}$ such that $s_0=ord_{d_x}(U)\leq s$. If $U=0$ we are done, so suppose that $U\neq 0$.
We have
\begin{equation}\label{3.1}
U=ZQ=T(d_x+ad_y+b_{n+1})
\end{equation}

\noindent
for suitable $Z,T\in F[d_x,d_y]$, hence $s_0=s$ and $ord_{d_x}(Z)=0, ord_{d_x}(T)=s-1$. One can expand $T=t_{s-1}E^{s-1}+\cdots +
t_0$ for appropriate $t_0,\dots, t_{s-1}\in F[d_y]$. Thus, the equation (\ref{3.1}) one rewrite with respect to the powers
of $E$:
$$Z(q_sE^s+\cdots +q_0)=(t_{s-1}E^{s-1}+\cdots +t_0)(E+b_{n+1})$$
\noindent
which is equivalent to a system of the following $s+1$ equalities:

\begin{equation}\label{3.2.j}
Zq_j=t_{j-1}+t_jb_{j,j}+t_{j+1}b_{j,j+1}+\cdots +t_{s-1}b_{j,s-1}
\end{equation}

\noindent
for suitable $b_{j,j},\dots, b_{j,s-1} \in F; 1\leq j\leq s$ and

\begin{equation}\label{3.2.0}
Zq_0=t_0b_{n+1}+t_1b_{0,1}+t_2b_{0,2}+\cdots +t_{s-1}b_{0,s-1}
\end{equation}

\noindent
Viewing the right-hand sides of the equations (\ref{3.2.j}), (\ref{3.2.0}) as a linear system in $t_0,\dots , t_{s-1}$ we
get that there is a unique linear combination (from the right) of $s$ expressions in the right-hand sides
of (\ref{3.2.j}) which equals (\ref{3.2.0}), the
coefficients $f_1,\dots, f_s$ of this combination belong to $F$. Therefore, the solvability of (\ref{3.2.j}), (\ref{3.2.0})
in $Z\neq 0,t_0,\dots, t_{s-1}$ entails the equality

\begin{equation}\label{3.3}
q_1f_1+\cdots +q_sf_s=q_0
\end{equation}

\noindent
Thus (\ref{3.1}) implies (\ref{3.3}). Hence as a solution of the system (\ref{3.2.j}),
(\ref{3.2.0}) one can take $Z=1$ and consecutively express $t_{s-1}\in F$ from the equation
(\ref{3.2.j}) with $j=s$, after that express $t_{s-2} \in F$ from the equation (\ref{3.2.j})
with $j=s-1$ and so on, finally express $t_0\in F$ from (\ref{3.2.j}) with $j=1$. The last
equation (\ref{3.2.0}) of the system is fulfilled due to (\ref{3.3}). As a result we obtain
(cf. (\ref{3.1})) $Q=(t_{s-1}E^{s-1}+\cdots +t_0)(E+b_{n+1})$ with $t_i\in F$, in other
words $ {\hat {I_1}}\cap \dots \cap {\hat I_n} = \langle Q \rangle \subset {\hat I_{n+1}}
\subset F[E]$.

This leads to contradiction with the assumption ${\hat {I_1}}\cap \dots \cap {\hat I_n} \supsetneq
{\hat {I_1}}\cap \dots \cap {\hat{I_n}} \cap
{\hat I_{n+1}}$, which shows that the supposition $U\neq 0$ was wrong, thus ${{I_1}}\cap \dots \cap {{I_n}} \cap
{I_{n+1}} = \langle M \rangle$. The proposition is proved. \bull

\begin{corollary}\label{principal}
The ideal ${ {I_1}}\cap \dots \cap {{I_n}}$ is generated by an element from $F[E]$.
\end{corollary}

Now let the ideals $I_i=\langle d_x+a_id_y+b_i \rangle \subset F[d_x,d_y]$ be given, where
$a_i,b_i\in F, 1\leq i\leq k$. Our goal is to study their intersection $I={ {I_1}}\cap \dots \cap {{I_k}}$. Combining together all the classes of the ideals with the same $a_i$ and making use
of Corollary~\ref{principal} we replace the intersection from one class by $\langle Z_i
\rangle$ for a certain $Z_i\in F[E_i]$ where $E_i=d_x+a_id_y$. Then $I={ {I_1}}\cap \dots \cap {{I_k}}=\langle Z_1 \rangle \cap \cdots \cap \langle Z_l \rangle$ for some $l$. Denote
$s_i=ord (Z_i); 1\leq i\leq l$ and $s=s_1+\cdots + s_l$.

\begin{lemma}\label{order}
For any $Q\in I$ we have $ord_{d_x} (Q) \geq s$.
\end{lemma}

{\bf Proof.} Observe that $symb(Q)$ is divided by
$\prod _{1\leq i\leq l} (d_xf + a_id_yf)^{s_i}$ treated as a homogeneous polynomial in
$d_xf, d_yf$. \bull

\begin{theorem}\label{criterium}
a) The ideal $I$ is principal if and only if $I$ contains $Q$ with the order $ord(Q)\leq s$;

b) in this case $ord(Q)=s$ and $I=\langle Q \rangle$.
\end{theorem}

{\bf Proof.} Obviously, the typical differential dimension $dim(\langle Z_i \rangle)=s_i;
1\leq i\leq l$ \cite{Kolchin} and $dim(I)\leq s$ due to \cite{Cassidy, Sit}. Hence if
$I=\langle L \rangle$ is principal then $ord(L)=dim(I)\leq s$.

Conversely, let $Q\in I$ and $ord(Q)\leq s$, by virtue of Lemma~\ref{order} we have $ord(Q)=s$
and the derivative $d_x^s$ occurs in $Q$. Our purpose is to show that $I=\langle Q \rangle$.
Indeed, take any $V\in I$ and divide $V$ by $Q$ with remainder, we get $V=WQ+U$ where
$ord_{d_x}(U)<s$, therefore, $U=0$ due to Lemma~\ref{order}. Thus, $I=\langle Q \rangle$. \bull

\begin{corollary}
Let the differential field $F=\QQ(x,y)$. There is a {\it polynomial-time} algorithm which tests
whether $I$ is principal.
\end{corollary}

{\bf Proof.} First the algorithm produces $Z_i; 1\leq i\leq l$ by finding a non-zero solution
of a linear (algebraic) homogeneous system on the coefficients from $F$ of $T_1,\dots , T_n\in
F[E_i]$ such that $T_1(d_x+a_id_y+b_1)=\cdots =T_n(d_x+a_id_y+b_n)$ with the minimal possible
order $ord(T_1)=\cdots =ord(T_n)$ (trying consecutively the orders 1,2...). Denote
$Z_i=T_1(d_x+a_id_y+b_1), s_i=ord(Z_i)$, then $Z_i$ is a generator of the ideal $\langle d_x+a_id_y+b_1
\rangle \cap \cdots \cap \langle d_x+a_id_y+b_n \rangle$, see Corollary~\ref{principal}.

Thereupon the algorithm looks for $V_1,\dots, V_l \in F[d_x,d_y]$ with $ord(V_i)\leq s-s_i;
1\leq i\leq l$ such that $V_1Z_1=\cdots = V_lZ_l$. The latter we treat as a linear (algebraic)
homogeneous system in the coefficients from $F$ of $V_1,\dots, V_l$. Theorem~\ref{criterium}
entails that this system has a non-zero solution if and only if $I$ is principal. \bull

\begin{remark}
Observe that the usual method of finding the intersection of ideals invoking Groebner bases,
runs in double-exponential time.
\end{remark}

\subsection{Constructing intersection of all first-order factors}\label{radical}

In this subsection $F$ denotes a universal field \cite{Kolchin} with two derivatives $d_x,d_y$.

The purpose of this subsection is to construct the intersection $U\subset F[d_x,d_y]$ of all the
principal ideals $\langle L \rangle$ for the first-order factors $L \in F[d_x,d_y]$ of $T\in \QQ (x,y)[d_x,d_y]$.
Evidently, $U\supset  \langle T \rangle$. We mention that in \cite{GS05} a radical of a module
of a differential type $\tau$ was defined as the intersection of the maximal classes of
$\tau$-equivalent modules, and a question was posed whether one can calculate the radical. Here
$U$ (which could be called a {\it first-order radical}) is defined as an ideal (rather than a class
of equivalent ideals) and moreover, we calculate $U$.

Observe that the construction from the Subsection~\ref{factors} represents the family $V$ of all
the solutions of the form $hG$ (and which correspond to first-order factors of $T$ due to
Lemma~\ref{correspondence}) as follows (we use the notations from Subsections~\ref{factors},
~\ref{ideals}).
We assume that $a$ is fixed, while $f$ just satisfies
the equality $d_xf+ad_yf=0$. The family $V$ is a union of subfamilies of the form $V_0$ where
$V_0$ is given by means of a Janet base

\begin{equation}\label{4.1}
\{\sum_{i_1,i_2} v_{i_1,i_2,l} d_x^{i_1} d_y^{i_2} h\} _l
\end{equation}

\noindent
for $h$ where $v_{i_1,i_2,l}\in R$ together with a system (\ref{1.2}) for $f$.

For each element $hG\in V$ consider the first-order LPDO $L_{hG}=d_x+ad_y+b_{hG}$ such that
$L_{hG}(hG)=0$ (see Lemma~\ref{correspondence}). We claim that one can extend
Proposition~\ref{intersection} from a finite to an infinite number of principal ideals and
conclude that the ideal $\cap _{hG\in V} \langle L_{hG} \rangle$ is principal and moreover, is
generated by a suitable element $Q=\sum_{0\leq i\leq s} q_iE^i \in F[E]$ (see
Corollary~\ref{principal}). Indeed, one add consecutively the ideals $I_1=\langle L_{h_1G_1}
\rangle, I_2=  \langle L_{h_2G_2} \rangle, \dots$ for $h_jG_j\in V$, while the intersection
${\hat I}_1 \cap \cdots \cap {\hat I}_{j-1} \cap {\hat I}_j \subsetneq
{\hat I}_1 \cap \cdots \cap {\hat I}_{j-1}$ decreases (cf. the proof of
Proposition~\ref{intersection}). Then ${\hat I}_1 \cap \cdots \cap {\hat I}_{j}=
\langle Q_j=\sum _{0\leq i\leq j} q_{i,j} E^i\rangle$ for appropriate $q_{i,j}\in F$ (cf.
the proof of
Proposition~\ref{intersection}). Hence $\langle T\rangle \subset I_1\cap \cdots \cap I_j =
\langle Q_j\rangle$ due to Corollary~\ref{principal}. Thus, $j\leq n$ and
$\cap _{hG\in V} \langle L_{hG} \rangle = I_1\cap \cdots \cap I_j$ which proves the claim.

To produce $Q=Q_j=\sum_{0\leq i\leq j} q_iE^i$ the algorithm successively tries $j=0,1,\dots$,
treating $q_i$ as indeterminates. The aim is to find $Q$ such that $Q(hG)=0$ for any $hG\in V_0$
(for each subfamily $V_0$ of $V$). The algorithm expands $Q(hG)=A_0G^{(0)}+\cdots +A_jG^{(j)}$ (cf.
Subsection~\ref{factors}). One can view each $A_i$ as an LPDO in $h$ with the coefficients being
linear forms in $q_0,\dots, q_j$ over $R$. The algorithm divides every $A_i, 0\leq i\leq j$ with
the remainder by the Janet base (\ref{4.1}), as a result we obtain LPDO ${\bar A_i}=\sum_{i_1,i_2}
a_{i,i_1,i_2}d_x^{i_1}d_y^{i_2}$. Thus, $Q$ vanishes at any $hG\in V_0$ if and only if $a_{i,i_1,i_2}=0$
for all $0\leq i\leq j; i_1,i_2$ under condition (\ref{1.2}).

Denote by $\cal S$ the conjunction of the systems $a_{i,i_1,i_2}=0$ for all $0\leq i\leq j; i_1,i_2$
and for all subfamilies of the form $V_0$ of $V$. One can treat $\cal S$ as a homogeneous linear over
$q_0,\dots, q_j$ system  with parameters being derivatives $f,d_yf,\dots,d_y^lf$ for a certain $l$.
Solving this parametric linear system (see e.g. \cite{G90}) the algorithm finds the (algebraic)
conditions on $f,d_yf,\dots,d_y^lf$ under which the system is solvable and in addition, finds the
expressions for solutions (being rational functions in the parameters). After that the algorithm tests
whether these conditions are compatible  with  (\ref{1.2}), applying the subroutine from the
elimination procedure which yields formula (\ref{1.2}) in Subsection~\ref{factors}. If yes
then the algorithm produces a solution
$q_0,\dots, q_j\in R$ of the parametric linear system. Else, the algorithm proceeds from the current
value $j$ to the next value $j+1$.

Thus, the algorithm for each $a$ such that $d_xf+ad_yf$ is a (linear) divisor of $symb(T)$, produces
applying the described above construction a generator $Q_a\in F[d_x,d_y]$ of the (principal) ideal
being the intersection of all the principal ideals generated by the divisors of the form
$d_x+ad_y+b$ of $T$ for varying $b$. Finally, the algorithm finds the intersection
$U=\cap _a \langle Q_a \rangle$ over all the divisors $d_xf+ad_yf$ of $symb(T)$ making use of Janet bases
(cf. \cite{GS05}). Thus, the following theorem is proved.

\begin{theorem}
For any LPDO $T\in \QQ(x,y)[d_x,d_y]$ one can construct the intersection of all the principal ideals
generated by the first-order factors of $T$.
\end{theorem}

\section{Fractional-derivatives series solutions of a second-order operator
and factoring}
\label{second}

In this section we study a particular case of a second-order LPDO
$T=T_0+T_1+T_2$ and describe its possible fractional-derivatives series
solutions being outputs of the algorithm from Section~\ref{solution}.
First, if the symbol $symb(T)$ is separable then for each of its two
different linear divisors $d_xf_1+ad_yf_1$ the algorithm provides a
fractional-derivatives series solution of $T$ of the form (cf. (\ref{2.1}))
$$\sum_{0\leq i< \infty} h_iG^{(-i)}$$
\noindent
where $G=G(f_1)$ and $d_xf_1+ad_yf_1=0, \quad f_1\neq const$.

From now on let us assume that $symb(T)$ is non-separable and write
$T=d_x^2+2ad_xd_y+a^2d_y^2+b_{0,1}d_x+b_{1,0}d_y+b_{0,0}$. The first
step of the algorithm from Section~\ref{solution} yields $f_1$ such that
$d_xf_1+ad_yf_1=0$. Introduce the discriminant of $T$ as follows:
$$(-T+b_{0,0})f_1=(d_xa+2ad_ya+ab_{0,1}-b_{1,0})(d_yf_1)=Disc\cdot (d_yf_1).$$

If $Disc\neq 0$ we take any $f_2$ which satisfies the following (non-homogeneous) first-order LPDE:
$$d_xf_2+ad_yf_2=\sqrt{Disc\cdot (d_yf_1)}$$
\noindent
and $G=G_{1/2}(f_1,f_2)$ (see Definition~\ref{symbol}), then the algorithm
constructs a
fractional-derivatives series
$$\sum_{0\leq i<\infty} h_i G^{(-i/2)}$$
\noindent
being a solution of $T=0$. Each of two values of the sign of the square root
provides a generic solution of the multiplicity 1 (see
Section~\ref{generic}).
It corresponds to the leading edge with the
endpoints $(0,2),(1,0)$ having the slope $1/2$ at the second step of the
algorithm.

When $Disc=0$ the algorithm yields a (fractional-derivatives series) solution
$hG(f_1)$ of $T=0$ for an arbitrary particular $h$ such that $T(h)=0$. It
corresponds to the leading edge with the endpoints $(0,2),(0,0)$ having the
slope $0$ at the second step of the algorithm and provides a generic solution
of the multiplicity $2$. Relying on Lemma~\ref{correspondence} one obtains
the following corollary (cf. \cite{GS04}).

\begin{corollary}
A second-order LPDO with a non-separable symbol is irreducible if and only if
$Disc \neq 0$.
\end{corollary}

\vspace*{0.4 cm} {\bf Acknowledgements.}
The author is grateful to the Max-Planck Institut fuer Mathematik, Bonn
where the paper was written,
to Y.Manin for his attention to the work and to F.Schwarz who has
pointed to the Laplace method.


\begin{thebibliography}{99}
\bibitem{Cano}
F.~Aroca, J.Cano, {\em Formal solutions of linear PDEs and convex polyhedra},
J. Symb. Comput., {\bf 32} (2001), 717-737.
\bibitem{Bjork}
J.-E.~Bj\"ork, {\em Rings of differential operators}, North-Holland, 1979.
\bibitem{Cassidy} P.~Cassidy, {\em Differential Algebraic Groups}, Amer. J.
 Math., {\bf 94}, (1972), 891--954
\bibitem{Goursat}
E.~Goursat, {\em Le\c con sur l'int\'egration des \'equations aux
d\'eriv\'ees partielles, vol. I, II}, A.Hermann, 1898.
\bibitem{G86}
D.~Grigoriev, {\em Computational complexity in polynomial algebra}, in
Proc.Intern.Congress Mathem., Berkeley (1986), 1452--1460.
\bibitem{G89}
D.~Grigoriev, {\em Complexity of quantifier elimination in the theory of ordinary differential
equations}, Lect. Notes Comput. Sci., {\bf 378}, (1989), 11--25.
\bibitem{G90}
D.~Grigoriev, {\em Complexity of factoring and calculating the GCD of linear ordinary
differential operators}, J. Symbolic Computations,
{\bf 7} (1990), 7--37.
\bibitem{G05}
D.~Grigoriev, {\em Weak B\'ezout inequality for $D$-modules},
J.Complexity, {\bf 21} (2005), 532--542.
\bibitem{GS04} D.~Grigoriev, F.~Schwarz, {\em Factoring and solving linear
partial differential equations}, Computing~{\bf 73} (2004), 179--197
\bibitem{GS05}
D.~Grigoriev, F.~Schwarz, {\em Loewy- and primary-decompositions of
$D$-modules}, Adv. Appl. Math. {\bf 38} (2007), 526--541.
\bibitem{GS91}
D.~Grigoriev, M.~Singer, {\em Solving ordinary differential
equations in the series with real exponents}, Trans. AMS, {\bf
327} (1991), 329--351.
\bibitem{Janet}
M.~Janet, {Les modules de formes alg\'ebriques et la th\'eorie
g\'en\'erale des syst\`emes diff\'erentiels}, Annals Sci. Ecole
Normale Sup\'er., {\bf 41} (1924), 27--65.
\bibitem{Kolchin}
E.~Kolchin, {\em Differential algebra and algebraic groups}, Academic Press,
1973.
\bibitem{Pankratiev}
M.~Kondratieva, A.~Levin, A.~Mikhalev, E.~Pankratiev, {\em Differential and
difference dimension polynomials}, Kluwer, 1999.
\bibitem{Lakshman}
Y.~Lakshman, D.~Saunders, {\em Sparse polynomial interpolation in
nonstandard bases}, SIAM J. Comput, {\bf 24} (1995), 387--397.
\bibitem{Sabbah}
{\em D-modules coh\'erents et holonomes}, Travaux en cours, {\bf 45},
eds. P.~Maisonobe, C.~Sabbah, Hermann, Paris, 1993.
\bibitem{Malgrange81}
B.~Malgrange, {\em R\'eduction d'un syst\`eme microdiff\'erentiel au point
g\'en\'erique}, Compositio Mathematica, {\bf 44} (1981), 133--143.
\bibitem{Malgrange91}
B.~Malgrange, {\em Equations diff\'erentielles \`a coefficients polynomiaux},
Progress in Math., {\bf 96}, Birkhauser, 1991.
\bibitem{Singer2003} M.~van der Put, M.~Singer, {\em Galois theory of linear
differential equations}, Grundlehren der Mathematischen Wissenschaften,
{\bf 328},
Springer, 2003.
\bibitem{Quadrat} A.~Quadrat, {\em An introduction to the algebraic theory of linear systems of partial differential equations}, www-sop.inria.fr/cafe/Alban.Quadrat/Temporaire.html
\bibitem{Sturmfels}
M.~Saito, B.~Sturmfels, N.~Takayama, {\em Gr\"obner deformations of hypergeometric
differential equations}, Algorithms and Computation in Mathematics, {\bf 6},
Springer, 2000.
\bibitem{schwarz}
F.~Schwarz, {\em Janet bases for symmetry groups, Groebner bases
and applications}, in London Math. Society, Lecture Notes Ser.,
{\bf 251}, Cambridge University Press, (1998), 221--234.
\bibitem{S56} A.~Seidenberg {\em An elimination theory for differential algebra}, Univ. Calif.
Publ. Math., (N.S.), {\bf 3}, (1956), 31--65
\bibitem{Sit}
W.Yu.~Sit, {\em Typical differential dimension of the intersection of linear
differential algebraic groups}, J.Algebra {\bf 32} (1974), 476--487.
\bibitem{Tsarev} S.~Tsarev, {\em Factorization of linear partial differential operators and the Darboux method for integrating nonlinear partial differential equations}, Theoret. and Math. Phys., {\bf 122}, (2000), 121--133.
\bibitem{Walker}
R.~Walker, {\em Algebraic curves}, Princeton, 1950.
\bibitem{Wasow}
W.~Wasow, {\em Asymptotic expansions for ordinary differential equations}, New York,
Krieger Publ. Co., 1976.
\end{thebibliography}
\end{document}